\documentclass[12pt,reqno]{amsart}

\usepackage{amscd}
\usepackage{amssymb}

\newcommand{\N}{{\mathbb N}}
\newcommand{\Z}{{\mathbb Z}}
\newcommand{\Q}{{\mathbb Q}}
\newcommand{\C}{{\mathbb C}}

\renewcommand{\P}{{\mathbb P}}
\renewcommand{\H}{{\mathbb H}}

\newcommand{\II}{{\mathcal I}}
\newcommand{\KK}{{\mathcal K}}

\newcommand{\NN}{{\mathcal N}}

\newcommand{\OO}{{\mathcal O}}
\newcommand{\PP}{{\mathcal P}}

\newcommand{\TT}{{\mathcal T}}

\newcommand{\ddd}{{\rm d}}
\newcommand{\www}{\widetilde}

\newcommand{\paa}{\partial}

\newcommand{\mmod}{{\rm mod}}

\DeclareMathOperator{\Aut}{Aut}
\DeclareMathOperator{\codim}{codim}

\DeclareMathOperator{\id}{id}

\DeclareMathOperator{\Lie}{Lie}
\DeclareMathOperator{\mult}{mult}

\theoremstyle{plain}
\newtheorem{lemma}{Lemma}[section]
\newtheorem{theorem}[lemma]{Theorem}

\newtheorem{corollary}[lemma]{Corollary}

\theoremstyle{definition}
\newtheorem{definition}[lemma]{Definition}

\newtheorem{remark}[lemma]{Remark}
\newtheorem{remarks}[lemma]{Remarks}

\newtheorem{examples}[lemma]{Examples}

\newtheorem{definition/lemma}[lemma]{Definition/Lemma}
\newtheorem{notations}[lemma]{Notations}

\begin{document}

\title[3-dimensional $F$-manifolds]
{3-dimensional $F$-manifolds} 

\author[A. Basalaev and C. Hertling]{Alexey Basalaev and Claus Hertling}
\address{Faculty of Mathematics, National Research University Higher School of Economics, Usacheva str., 6, 119048 Moscow, Russian Federation, and Skolkovo Institute of Science and Technology, Nobelya str., 3, 121205 Moscow, Russian Federation}
\email{a.basalaev\char64 skoltech.ru}
\address{Lehrstuhl f\"ur algebraische Geometrie,
Universit\"at Mannheim, B6, 26, 68159 Mannheim, Germany}
\email{hertling\char64 math.uni-mannheim.de}

\dedicatory{In remembrance of Boris Dubrovin}

\thanks{This work was funded by the Deutsche Forschungsgemeinschaft
(DFG, German Research Foundation) -- 242588615}



\date{December 20, 2020}

\begin{abstract}
$F$-manifolds are complex manifolds with a multiplication
with unit on the holomorphic tangent bundle with a certain 
integrability condition. Here the local classification
of 3-dimensional $F$-manifolds with or without Euler fields
is pursued. 
\end{abstract}

\maketitle


\tableofcontents


\section{Introduction}\label{c1}
\setcounter{equation}{0}

\noindent 
Boris Dubrovin defined and studied 
{\it Frobenius manifolds} \cite{Du92} \cite{Du96}.
A Frobenius manifold is a complex manifold $M$ with
a holomorphic flat metric $g$ and a holomorphic commutative
and associative multiplication $\circ$ with unit $e$ on the
holomorphic tangent bundle $TM$ such that 
$g(X\circ Y,Z)=g(X,Y\circ Z)$ and such that locally 
a holomorphic function $\Phi$ (a {\it potential})
with $g(X\circ Y,Z)=XYZ(\Phi)$ for flat vector fields
$X,Y,Z$ exists. Often one has additionally an Euler field
$E$, a holomorphic vector field with 
$\textup{Lie}_E(\circ)=1\cdot \circ$ and 
$\textup{Lie}_E(g)=D\cdot g$ for some $D\in\C$. 

This seemingly purely differential geometric object
has many different facets and lies at the crossroads
of very different mathematical areas, integrable systems,
meromorphic connections, singularity theory, quantum cohomology
and thus mirror symmetry.
Boris Dubrovin explored many of these crossroads.

Manin and the second author defined the notion of
an {\it $F$-manifold} \cite{HM99}. 
It is a complex manifold $M$ with a holomorphic commutative
and associative multiplication $\circ$ with a unit $e$
on the holomorphic tangent bundle which satisfies the
integrability condition 
\begin{eqnarray}\label{1.1}
\Lie_{X\circ Y}(\circ)= X\circ\Lie_Y(\circ)+Y\circ\Lie_X(\circ)
\textup{ for }X,Y\in\OO(TM).\hspace*{0.5cm}
\end{eqnarray}
Here an Euler field is a holomorphic vector field $E$
with $\textup{Lie}_E(\circ)=1\cdot \circ$. 

Frobenius manifolds are $F$-manifolds, and this is the
original motivation for the definition of $F$-manifolds. 
But there are also $F$-manifolds which cannot be enriched 
to Frobenius manifolds. 
The paper \cite{DH20-2} starts with $F$-manifolds and
studies how and when they can be enriched to Frobenius
manifolds. Crucial is the existence of a certain bundle with
a meromorphic connection (called $(TE)$-structure in 
\cite{DH20-2}) over an $F$-manifold.

Slightly weaker, but almost as strong as a Frobenius manifold
is the notion of a {\it flat $F$-manifold}, which was defined
by Manin \cite{Ma05}. It is an $F$-manifold with flat 
connection $D$ on $TM$ with $D(C^M)$ and $D(e)=0$, where
$C^M$ is the Higgs field from the multiplication, so 
$C^M_X=X\circ:TM\to TM$ for $X\in\OO(TM)$.
Then an Euler field $E$ is an Euler field of the $F$-manifold
such that $D_\bullet E:TM\to TM$ is a flat endomorphism.

Recently, flat $F$-manifolds with Euler fields were subject 
to work by Arsie and Lorenzoni \cite{AL13} \cite{Lo14}
\cite{AL17} \cite{AL19}, 
Kato, Mano and Sekiguchi \cite{KMS15}, Kawakami and Mano
\cite{KM19}, Konishi, Minabe and Shiraishi \cite{KMS18}
\cite{KM20}. They established such structures on orbit spaces
of complex reflection groups. 
And especially they observed a beautiful correspondence between
{\it regular} flat 3-dimensional $F$-manifolds
and solutions of the Painlev\'e equations of types
VI, V and IV \cite{AL19} \cite{KMS15} \cite{KM19}.

A {\it regular} $F$-manifold is an $F$-manifold with
Euler field such that the endomorphism $E\circ$ on $TM$
has everywhere for each eigenvalue only one Jordan block.
This notion was defined and studied by David and the
second author \cite{DH17}. 

The second author studied $F$-manifolds in
\cite[ch. 1--5]{He02}. There he classified all germs of
2-dimensional $F$-manifolds with or without Euler fields.
This classification is easy, see below. 
But already the classification of the germs of 3-dimensional
$F$-manifolds is rich. It was not pursued systematically
in \cite{He02} or anywhere else. 

This paper aims at a systematic classification of
germs of 3-dimensional $F$-manifolds. It succeeds in the majority of the cases, but not in all cases. 
One motivation is the interest in regular flat 3-dimensional
$F$-manifolds. The classification in this paper gives 
especially all germs of regular 3-dimensional $F$-manifolds.

In order to distinguish different cases,
the 3-dimensional algebras over $\C$ have to be listed.

\begin{remarks}\label{t1.1}
Here the commutative and associative algebras with unit over
$\C$ of dimensions 1, 2 and 3 are listed.
In dimension 1, the only algebra is $\C$.
In dimension 2, there exist up to isomorphism two algebras 
\begin{eqnarray*}
P^{(1)}&:=&\C[x]/(x^2),\\
P^{(2)}&:=&\C\oplus \C.
\end{eqnarray*}
In dimension 3, there exist up to isomorphism four algebras,
\begin{eqnarray*}
Q^{(1)}&:=& \C[x_1,x_2]/(x_1^2,x_1x_2,x_2^2),\\
Q^{(2)}&:=& \C[x]/(x^3),\\
Q^{(3)}&:=& \C\oplus \C[x]/(x^2)=\C\oplus P^{(1)},\\
Q^{(4)}&:=& \C\oplus \C\oplus \C
\end{eqnarray*}
A sum $\bigoplus_{j=1}^n\C$ of 1-dimensional algebras is called 
{\it semisimple}, so $\C$, $P^{(2)}$ and $Q^{(4)}$ are 
semisimple. 
The algebras $\C$, $P^{(1)}$, $Q^{(1)}$ and $Q^{(2)}$ 
are irreducible. The decomposition of each algebra into
irreducible algebras is unique. 
The algebras $\C$, $P^{(1)}$, $P^{(2)}$, $Q^{(2)}$, 
$Q^{(3)}$ and $Q^{(4)}$ are Gorenstein rings, $Q^{(1)}$ 
is not a Gorenstein ring. 
(The notations $Q^{(1)}$ and $Q^{(2)}$ are opposite
to those in \cite[5.5]{He02}.)
\end{remarks}

Now let $(M,\circ,e)$ be a connected 3-dimensional complex 
manifold with a commutative and associative multiplication on 
$TM$ with unit field $e$, but not necessarily with \eqref{1.1}.
Choose local coordinates $t=(t_1,t_2,t_3)$, and denote by
$y=(y_1,y_2,y_3)$ the fiber coordinates on $T^*M$ such that
$y_j$ corresponds to the coordinate vector field 
$\paa_j:=\paa/\paa t_j$. Then $\alpha:=\sum_{j=1}^3y_j\ddd t_j$
is the canonical 1-form on $T^*M$. The multiplication
is given by $\paa_i\circ\paa_j=\sum_{k=1}^3a_{ij}^{k}\paa_k$
with coefficients $a_{ij}^k\in\OO_M$. It gives rise to
the sheaf of ideals $\II_M\subset\OO(T^*M)$ with
\begin{eqnarray}\label{1.2}
\II_M:=\bigl(y_iy_j-\sum_{k=1}^3a_{ij}^{k}y_k
\,|\,i,j\in\{1,2,3\}\bigr)\subset\OO(T^*M)
\end{eqnarray}
and the complex space $L_M\subset T^*M$ which is as a set
the zero set of $\II_M$ and which has the complex structure
$\OO_{L_M}=(\OO_{T^*M}/\II_M)|_{L_M}$. 
The projection $\pi_L:L_M\to M$ is
flat and finite of degree 3. For each $t\in M$, the 
points in $\pi_L^{-1}(t)\subset T_t^*M$ are the 
simultaneous eigenvalues of all endomorphisms 
$X|_t\circ:T_tM\to T_tM$ for $X\in T_tM$. They 
correspond to the irreducible subalgebras of $T_tM$. 

The numbering $Q^{(1)},...,Q^{(4)}$ above was chosen so that
for each $j\in\{1,2,3,4\}$, the subset
$\bigcup_{i\leq j}\{t\in M\,|\,T_tM\cong Q^{(i)}\}$ is 
empty or an analytic subvariety of $M$ or equal to $M$
(Lemma \ref{t4.3} gives more precise statements). 
The algebra $Q^{(j)}$ with $T_tM\cong Q^{(j)}$ for generic 
$t\in M$ is called the {\it generic type} of $M$
($M$ is connected). $M$ is called {\it generically semisimple}
if the generic type is $Q^{(4)}$. 

$L_M$ is called {\it analytic spectrum} of $(M,\circ,e)$.
It encodes the multiplication and is crucial for its
understanding. The integrability condition \eqref{1.1}
of an $F$-manifold is equivalent to 
$\{\II_M,\II_M\}\subset\II_M$, where $\{.,.\}$ is the 
Poisson bracket on $\OO(T^*M)$ \cite{HMT09} (cited in
Theorem \ref{t2.12}). In the generically semisimple
case, this is equivalent to $L_M^{reg}\subset T^*M$ being
Lagrange. This connects the generically semisimple
$F$-manifolds with the Lagrange fibrations and Lagrange 
maps of Arnold \cite[ch. 18]{AGV85}. Givental's paper 
\cite{Gi88} on Lagrange maps contains implicitly many results
and examples of generically semisimple $F$-manifolds. 

In the case of an $F$-manifold, the integrability condition
\eqref{1.1} implies that at a point  $t\in M$ such that
$T_tM$ decomposes into several irreducible algebras,
also the germ of the $F$-manifold decomposes uniquely
into a product of germs of $F$-manifolds, one for each
summand of $T_tM$ \cite[Theorem 2.11]{He02} (cited in 
Theorem \ref{t2.5}), and an Euler field decomposes
accordingly. Therefore in the classification of germs of
$F$-manifolds, we can restrict to the classification of
the irreducible germs, which are the germs $(M,0)$
such that $T_0M$ is irreducible. 
A rough distinction of classes is given by the isomorphism class 
of $T_0M$ and the generic type. 

The following table shows which examples, lemmas and theorems
in this paper concern which class of irreducible germs $(M,0)$ 
of $F$-manifolds of dimensions 1 or 2 or 3.

\begin{figure}[ht]
\begin{tabular}{c|c|l}
$T_0M$ & generic type & \\ \hline 
$\C$ & $\C$ & only 1 $F$-manifold $A_1$: Lemma \ref{t2.6}\\ 
\hline 
$P^{(1)}$, & $P^{(1)}$ & 1 $F$-manifold $\NN_2$: 
Theorems \ref{t3.1}, \ref{t3.2}\\ 
$P^{(1)}$ & $P^{(2)}$ & 1 series $I_2(m)$, $m\geq 3$:
Theorem \ref{t3.1} \\ \hline 
$Q^{(1)}$ & $Q^{(1)}$ & Theorem \ref{t5.2} \\
$Q^{(1)}$ & $Q^{(2)}$ & Theorem \ref{t5.4} b)+(c)\\
$Q^{(2)}$ & $Q^{(2)}$ & Theorem \ref{t5.4} (a) \\
$Q^{(1)}$ & $Q^{(3)}$ & ??, Lemma \ref{t5.8} \\
$Q^{(2)}$ & $Q^{(3)}$ & Theorem \ref{t5.6} \\
$Q^{(1)}$ & $Q^{(4)}$ & ??, Theorem \ref{t6.3}, 
Lemmas \ref{t6.4}, \ref{6.5}\\ 
$Q^{(2)}$ & $Q^{(4)}$ & Examples \ref{t6.2}, 
Theorems \ref{t6.3}, \ref{t7.1}
\end{tabular}
\caption[Table 1]{Table of results}
\end{figure}

The results for dimension 1 and 2 are cited from \cite{He02},
and they are easy. The classification in dimension 3 is 
surprisingly rich. The cases with $T_0M\cong Q^{(2)}$
are easier than those with $T_0M\cong Q^{(1)}$.
In the two cases with $T_0M\cong Q^{(1)}$ and 
generic type $Q^{(3)}$ or $Q^{(4)}$, we have no complete
classification, but just some examples. 
But as $Q^{(1)}$ does not turn up in a regular $F$-manifold, 
our classification gives all germs of regular $F$-manifolds.
Most not generically semisimple $F$-manifolds appear here
for the first time. And also most of the generically 
semisimple $F$-manifolds, namely those in
Theorem \ref{t7.1} with $T_0M\cong Q^{(2)}$, are new. 
Their classification is linked to the classification 
of germs of plane curves of multiplicity 3 
(see the Remarks \ref{t7.3}). 

Possibly the most interesting germs $(M,0)$ of 
3-dimensional $F$-manifolds 
are the generically semisimple germs with Euler field.
Those with $T_0M\cong Q^{(2)}$ are given in
Corollary \ref{t7.2}. 

Section \ref{c2} collects general facts on $F$-manifolds
from \cite{He02}. Section \ref{c3} recalls the classification
of the 2-dimensional germs of $F$-manifolds. 
Section \ref{c4} provides basic formulas for 
3-dimensional $F$-manifolds. 
Section \ref{c5} classifies the not generically semisimple
germs (except those with $T_0M\cong Q^{(1)}$ and generic
type $Q^{(3)}$). Section \ref{c6} gives examples
of generically semisimple $F$-manifolds.
Section \ref{c7} classifies the generically semisimple
germs with $T_0M\cong Q^{(2)}$.

\section{General facts on $F$-manifolds}\label{c2}
\setcounter{equation}{0}

\noindent 
$F$-manifolds were first defined in \cite{HM99}.
Their basic properties were developed in \cite{He02}.
This section reviews the main basic properties from 
\cite{He02} and an  additional fact from \cite{HMT09}.

\begin{definition}\label{t2.1} \cite{HM99}
(a) An {\it $F$-manifold} $(M,\circ,e)$ (without Euler field) 
is a holomorphic manifold $M$ with a holomorphic 
commutative and associative multiplication $\circ$
on the holomorphic tangent bundle $TM$ and with a 
global holomorphic vector field $e\in\TT_M:=\OO(TM)$ with
$e\circ=\id$ ($e$ is called a {\it unit field}),
which satisfies the integrability condition \eqref{1.1}.

(b) Given an $F$-manifold $(M,\circ,e)$, an {\it Euler field} on it is a global
vector field $E\in\TT_M$ with $\Lie_E(\circ)=\circ$.
\end{definition}

\begin{remark}\label{t2.2}
The integrability condition \eqref{1.1}
looks surprising at first sight. Though it is natural
from several points of view. Here are four of them.
\begin{list}{}{}
\item[(i)] Theorem \ref{t2.12} below rewrites condition
\eqref{1.1} as a natural condition on the ideal 
giving the analytic spectrum in $T^*M$.
\item[(ii)] Theorem \ref{t2.5} below gives a decomposition 
result for germs of $F$-manifolds. Condition \eqref{t2.1} 
is crucial in its proof in \cite{He02}.
\item[(iii)] The potentiality condition in a Frobenius 
manifold with holomorphic metric $g$ 
is equivalent to \eqref{1.1} plus the closedness of the
1-form (called coidentity) $g(e,.)$ \cite[Theorem 2.15]{He02}. 
\item[(iv)]
If the Higgs field of a $(TE)$-structure over a manifold $M$
is primitive, it induces on $M$ the structure of an
$F$-manifold with Euler field, see e.g. \cite{DH20-2}.
\end{list}
\end{remark}

\begin{remark}\label{t2.3}\cite[Proposition 2.10]{He02}
If one has $l$ $F$-manifolds $(M_k,\circ_k,e_k)$, 
$k\in\{1,...,l\}$, 
their product $M=\prod_{k=1}^l M_k$
inherits a natural structure of an $F$-manifold
$(M,\bigoplus_{k=1}^l \circ_k,
\sum_{k=1}^l (\textup{lift of }e_k\textup{ to }M)$. 
And if there are Euler fields $E_k$, then the sum 
$E=\sum_{k=1}^l(\textup{lift of }E_k\textup{ to }M)$
is an Euler field on the product $M$.
\end{remark}

\begin{remark}\label{t2.4}
A finite dimensional commutative and associative $\C$-algebra
$A$ with unit $e\in A$ decomposes uniquely into a direct sum
$A=\bigoplus_{k=1}^l A_k$ of local and irreducible
algebras $A_k$ with units $e_k$ with $e=\sum_{k=1}^l e_k$
and $A_{k_1}\circ A_{k_2}=0$ for $k_1\neq k_2$.
This is elementary (linear) algebra. The
decomposition is obtained as the simultaneous decomposition
into generalized eigenspaces of all endomorphisms
$a\circ:A\to A$ for $a\in A$ 
(see e.g. Lemma 2.1 in \cite{He02}).
The algebra $A$ is called {\it semisimple} if
$l=\dim A$ (so then $A_k=\C\cdot e_k$ for all $k$).
\end{remark}

Thanks to the condition \eqref{1.1}, this pointwise
decomposition extends in the case of an $F$-manifold
to a local decomposition, see Theorem \ref{t2.5}. 
This is the first important step in the local classification 
of $F$-manifolds.

\begin{theorem}\label{t2.5}\cite[Theorem 2.11]{He02}
Let $((M,t^0),\circ,e)$ be the germ at $t^0$ of an
$F$-manifold. 

\medskip
(a) The decomposition of the algebra 
$(T_{t^0}M,\circ|_{t^0},e|_{t^0})$ with unit into local
algebras extends into a canonical decomposition 
$(M,t^0)=\prod_{k=1}^l (M_k,t^{0,k})$ as a product
of germs of $F$-manifolds.

\medskip
(b) If $E$ is an Euler field of $M$, then $E$
decomposes as $E=\sum_{k=1}^l E_k$ with $E_k$
(the canonical lift of) an Euler field on $M_k$. 
\end{theorem}

\begin{lemma}\label{t2.6}\cite[Example 2.12 (i)]{He02}
In dimension 1, (up to isomorphism) there is only
one germ of an $F$-manifold, the germ $(M,0)=(\C,0)$ with
$e=\paa/\paa_{u_1}$, where $u_1$ is the coordinate on $\C$. 
Any Euler field on it has the shape 
$E=(u_1+c_1)e$ for some $c_1\in\C$. 
\end{lemma}

\begin{definition}\label{t2.7}
(a) Fix $n\in\N=\{1,2,..\}$ and define the set of its partitions,
\begin{eqnarray*}
\PP_n:=\{\beta=(\beta_1,...,\beta_{l(\beta)})\,|\,
\beta_i\in\N,\ \beta_i\geq\beta_{i+1},\sum_{i=1}^{l(\beta)}
\beta_i=n\}.
\end{eqnarray*}
For $\beta,\gamma\in\PP_n$ define
\begin{eqnarray*}
&&\beta\geq \gamma:\iff\exists\ \sigma:\{1,...,l(\gamma)\}
\to\{1,...,l(\beta)\}\textup{ s.t. }
\beta_j=\sum_{i\in\sigma^{-1}(j)}\gamma_i,\\
&&\beta>\gamma:\iff \beta\geq \gamma\textup{ and }\beta\neq
\gamma.
\end{eqnarray*}

(b) Let $(M,\circ,e)$ be an $F$-manifold of dimension $n$. 
Consider the map 
\begin{eqnarray*}
P:M\to\PP_n,\ P(t):=\{\textup{the partition of }n
\textup{ by the dimensions}\\
\textup{of the irreducible subalgebras of }T_tM\}\nonumber
\end{eqnarray*}

(c) An $F$-manifold is called {\it generically semisimple} 
if $P(t)=(1,...,1)$ ($\Leftrightarrow\beta(P(t))=n$) 
for generic $t$ 
(In \cite{He02} such an $F$-manifold is called {\it massive}). 
An $F$-manifold is called {\it semisimple} if it is semisimple
at all points.
\end{definition}

\begin{lemma}\label{t2.8}
Let $(M,\circ,e)$ be an $F$-manifold of dimension $n$.

(a) \cite[Proposition 2.5]{He02} 
For any $\beta\in\PP_n$, the set 
$\{t\in M\,|\, P(n)\geq \beta\}$ is an analytic subset of
$M$ or empty.

\medskip
(b) \cite[Proposition 2.6]{He02} 
Suppose that $M$ is connected. Then there is a unique
partition $\beta_0\in\PP_n$ such that the set
$\{t\in M\,|\, P(t)=\beta_0\}$ is open.
Its complement is called {\sf caustic} and is denoted
by $\KK:=\{t\in M\,|\, P(t)\neq\beta_0\}$. 
The caustic is an analytic hypersurface or empty. 
If $t\in\KK$, then $P(t)>\beta_0$. 

\medskip
(c) By Theorem \ref{t2.5} and Lemma \ref{t2.6}, 
a semisimple germ of an $F$-manifold is isomorphic to 
$(\C^n,0)$ with coordinates $u=(u_1,...,u_n)$ and partial units
$e_k=\paa_{u_k}$, which determine the multiplication by
$e_k\circ e_k=e_k$ and $e_{k_1}\circ e_{k_2}=0$
for $k_1\neq k_2$. The global unit field is $e=\sum_{k=1}^ne_k$.
The semisimple germ of dimension $n$ is said to be of type 
$A_1^n$. The coordinates $u_k$ or their shifts $u_k+c_k$
for any constants $c_1,..,c_n\in\C$ are 
{\sf Dubrovin's canonical coordinates}.
Any Euler field on this $F$-manifold has the shape 
$E=\sum_{k=1}^n (u_k+c_k)e_k$ for some $c_1,...,c_n\in\C$.
If an Euler field $E$ is fixed, the eigenvalues $u_k+c_k$ 
of $E\circ$ can be used as canonical coordinates.
This fixes their ambiguity.  
\end{lemma}

A generically semisimple $F$-manifold $M$ is easy locally
on $M-\KK$, but interesting near $\KK$. Then the 
canonical coordinates are coordinates locally on $M-\KK$.

Three notions from the theory of isolated hypersurface
singularities generalize to $F$-manifolds,
the {\it $\mu$-constant stratum}, the {\it modality}, 
and {\it simpleness}. 

\begin{definition}\label{t2.9}
Let $(M,\circ,e)$ be an $F$-manifold. 

\medskip
(a) For $p\in M$, the {\it $\mu$-constant stratum} of $p$
is the subvariety $S_\mu(p):=\{t\in M\,|\,
P(t)\geq P(p)\}.$ The {\it modality} $\mmod_\mu(M,p)$ is
\begin{eqnarray}\label{2.1}
\mmod_\mu(M,p):=\dim (S_\mu(p),p)-l(P(p)).
\end{eqnarray}

(b) The $F$-manifold is {\it simple} if $\mmod_\mu(M,p)=0$
for any $p\in P$. A simple $F$-manifold is generically 
semisimple because for any $F$-manifold 
$\mmod_\mu(p)=n-l(\beta_0)$ for $p\in M-\KK$. 
\end{definition}

The definition of the modality is motivated by the following.
If $(M,p)=\prod_{j=1}^{l(P(p))}(M_k,p^{(k)})$ as a germ of an 
$F$-manifold with idempotent vector fields
$e_1,...,e_{l(P(p))}$, then $\Lie_{e_j}(\circ)=0\cdot\circ$,
so the germs $(M,q)$ for $q$ in one integral manifold
of $e_1,...,e_{l(P(p))}$ are isomorphic as germs of 
$F$-manifolds. 

In the case of a generically semisimple $F$-manifold
with Euler field, the Euler field gives rise to a complementary
result.

\begin{theorem}\label{t2.10}\cite[Corollary 4.16]{He02}
Let $(M,\circ,e,E)$ be a generically semisimple $F$-manifold
with Euler field. For any $p\in M$, the set
\begin{eqnarray*}
\{t\in M\,|\, ((M,t),\circ,e,E)\cong ((M,p),\circ,e,E)\}
\end{eqnarray*}
is discrete and closed in $M$.
\end{theorem}

All the information of an $F$-manifold is carried also
by its analytic spectrum, which will be introduced now.

\begin{definition}\label{t2.11}
Let $(M,\circ,e)$ be a complex manifold of dimension $n$ 
with a holomorphic commutative and associative 
multiplication $\circ$ on the holomorphic 
tangent bundle and with a unit field $e$.
(In the moment, the condition \eqref{1.1}
is not imposed.)

\medskip
(a) We need some standard data on $T^*M$:
Let $\pi:T^*M\to M$ denote the projection.
Let $t=(t_1,...,t_n)$ be local coordinates
on $M$, and define $\paa_k:=\paa/\paa t_k$. 
Let $y=(y_1,...,y_n)$ be the fiber coordinates
on $T^*M$ which correspond to $(\paa_1,...,\paa_n)$.
Then the canonical 1-form
$\alpha$ takes the shape $\sum_{i=1}^ny_i{\rm d}t_i$,
and $\omega={\rm d}\alpha$ is the standard symplectic form.
The Hamilton vector field of 
$f\in\OO_{T^*M}$ is 
\begin{eqnarray}\label{2.2}
H_f = \sum_{k=1}^n \Bigl( \frac{\paa f}{\paa t_k}
\cdot \frac{\paa}{\paa y_k} 
-\frac{\paa f}{\paa y_k}\cdot \frac{\paa}{\paa t_k}\Bigr).
\end{eqnarray}
The Poisson bracket $\{.,.\}$ on $\OO_{T^*M}$ is defined by
\begin{eqnarray}\label{2.3}
\{f,g\}&:=& H_f(g)=\omega(H_f,H_g)=-H_g(f).
\end{eqnarray}

\medskip
(b) Define an ideal sheaf $\II_M\subset \OO_{T^*M}$
as follows. We choose coordinates $t_k$ and $y_k$ as in 
part (a) and such that $e_1=\paa_1$. 
Write
\begin{eqnarray}\label{2.4}
\paa_i\circ\paa_j =\sum_{k=1}^n a_{ij}^k\paa_k
\textup{ with }a_{ij}^k\in\OO_M.
\end{eqnarray}
Then 
\begin{eqnarray}\label{2.5}
\II_M:=\bigl(y_1-1,y_iy_j-\sum_{k=1}^na_{ij}^ky_k\bigr)
\subset\OO_{T^*M}.
\end{eqnarray}
The {\it analytic spectrum} (or {\it spectral cover})
$L_M:=\textup{Specan}_{\OO_M}(TM,\circ)\subset T^*M$
of $(M,\circ,e)$ 
is as a set the set at which the functions in $\II_M$ vanish.
It is a complex subspace of $T^*M$ with
complex structure given by 
$\OO_{L_M}=(\OO_{T^*M}/\II_M)|_{L_M}$.
\end{definition}

The analytic spectrum $L_M$ was studied in 
\cite[2.2 and 3.2]{He02}. 
But the following result was missed there.

\begin{theorem}\label{t2.12}\cite[2.5 Theorem]{HMT09}
A manifold $(M,\circ,e)$ with holomorphic commutative 
and associative multiplication $\circ$ on the 
holomorphic tangent bundle and unit field $e$ is an 
$F$-manifold if and only if $\{\II_M,\II_M\}\subset\II_M$.
\end{theorem}

\begin{remarks}\label{t2.13}
(i) The points in $L_M$ above a point $t\in M$ are
the 1-forms, which are the simultaneous eigenvalues
for all multiplication endomorphisms in $T_tM$. 
They are in 1-1 correspondence with the irreducible
subalgebras of $T_tM$.

\medskip
(ii) Let $(M,\circ,e)$ be a complex manifold of dimension
$n$ with commutative and associative 
multiplication on the holomorphic tangent bundle.
The projection $\pi|_{L_M}:L_M\to M$ is finite
and flat of degree $n$. The map 
\begin{eqnarray}\label{2.6}
{\bf a}:\TT_M\to\pi_*\OO(L_M),\quad X\mapsto 
\alpha(X)|_{L_M},
\end{eqnarray}
is an isomorphism of $\OO_M$-algebras.
In this way, the multiplication on $\pi_*\OO(L_M)$ determines
the multiplication on the tangent bundle.
The value $\alpha(X)(y,t)\in\C$ at a point $(y,t)\in L_M$
is the eigenvalue of $X\circ$ on the irreducible subalgebra
of $T_tM$ which corresponds to $(y,t)$. 

\medskip
(iii) In the case of a manifold with a multiplication
and unit field, such that the multiplication is  
generically semisimple, 
the restriction $L_M|_{M-\KK}$ of $L_M$ to $M-\KK$ 
is obviously smooth with $\dim M$ sheets above $M-\KK$.
Theorem 3.2 in \cite{He02} says that then $L_M$ is
reduced everywhere, so also above 
$L_M\cap \pi|_{L_M}^{-1}(\KK)$.

\medskip
(iv) In this situation, $\{\II_M,\II_M\}\subset\II_M$
says that $L_M$ is at smooth points a Lagrange submanifold
of $T^*M$. 

\medskip
(v) But in the case of a manifold with multiplication
and unit field, such that the multiplication is nowhere
semisimple, the analytic spectrum $L_M$ is nowhere
reduced. Then $\II_M$ is quite different from the
reduced ideal $\sqrt{\II_M}$.
Especially, the conditions
\begin{eqnarray}\label{2.7}
\{\II_M,\II_M\}\subset\II_M
\quad\textup{and}\quad  
\{\sqrt{\II_M},\sqrt{\II_M}\}\subset\sqrt{\II_M}
\end{eqnarray}
do not imply one another. 
The second condition in \eqref{2.7} is equivalent
to the condition that $L^{red}_M$ (the reduced space 
underlying $L_M$) is at smooth points a Lagrange submanifold
of $T^*M$. 
The examples 2.5.2 and 2.5.3 in \cite{HMT09}
and the examples below in Theorem \ref{t5.2} and Remark 
\ref{t5.3} (ii) with $b_2\neq 0$ 
are examples of $F$-manifolds (so 
$\{\II_M,\II_M\}\subset\II_M$ holds) with
$\{\sqrt{\II_M},\sqrt{\II_M}\}\not\subset\sqrt{\II_M}$.
The example (with $n=4$) in \cite[2.13 (v)]{DH20-2}
is an example of a manifold $(M,\circ,e)$ with 
$\{\II_M,\II_M\}\not\subset\II_M$ and 
$\{\sqrt{\II_M},\sqrt{\II_M}\}\subset\sqrt{\II_M}$.
\end{remarks}

We are mainly interested in the case of generically
semisimple $F$-manifolds. There the following result 
of Givental is relevant.

\begin{theorem}\label{t2.14}\cite[ch. 1.1]{Gi88}
An $n$-dimensional germ $(L,0)$ of a Lagrange variety
with embedding dimension $\textup{embdim}(L,0)=n+k$
with $k<n$ is a product of a $k$-dimensional
Lagrange germ $(L',0)$ with $\textup{embdim}(L',0)=2k$
and a smooth $(n-k)$-dimensional Lagrange germ
$(L'',0)$; here the decomposition of $(L,0)$
corresponds to a decomposition 
\begin{eqnarray*}
((S,0),\omega)\cong ((S',0),\omega')\times ((S'',0),\omega'')
\end{eqnarray*}
of the symplectic space germ $(S,0)$ which contains $(L,0)$.
\end{theorem}

Existence of an Euler field for a given $F$-manifold
is a problem with many facets. Some $F$-manifolds
have many Euler fields, others few, others none.
The cases of all 2- and many 3-dimensional $F$-manifolds
will be discussed in the sections \ref{c3} to \ref{c7}.
We are mainly interested in the generically semisimple
$F$-manifolds. There the following holds.

\begin{theorem}\label{t2.15}
(a) \cite[Theorem 3.3]{He02}
Let $(M,\circ,e)$ be a generically semisimple $F$-manifold.
A vector field $E$ is an Euler field if and only if
\begin{eqnarray}\label{2.8}
d({\bf a}(E))|_{L^{reg}_M}=\alpha|_{L_M^{reg}}.
\end{eqnarray}

(b) \cite[Lemma 3.4]{He02} 
Let $M$ be a sufficiently small representative of an
irreducible germ $(M,t^0)$ of a generically semisimple
$F$-manifold. For any $c\in\C$, there is a unique 
function $F:(L,(y^0,t^0))\to (\C,c)$ which is holomorphically
on $L^{reg}_M$ and continuous on $L$ (with value $c$
at $(y^0,t^0)$) and which satisfies 
$dF|_{L^{reg}_M}=\alpha|_{L^{reg}_M}$. 

\medskip
(c) The parts (a) and (b) imply that in the situation of (b),
for any $c\in \C$, there is a unique Euler field
$E_c$ on $M-\KK$ such that for $t\to t^0$
all eigenvalues of $E\circ$ tend to $c$. We have
$E_c=E_0+c\cdot e$. The characteristic polynomial
of $E_c\circ$ extends holomorphically to $t^0$
and has there the value $(x-c)^n$. 
The Euler field $E_c$ extends holomorphically to $M$
if and only if the function $F$ in part (b) is holomorphic
on $L_M$.
\end{theorem}

The question whether the function $F$ in part (b) is
holomorphic on $L_M$, will be rephrased in Theorem \ref{t2.16} 
(d), and a special case will be singled out in
Theorem \ref{t2.16} (e). 
Now we consider the germ $(S,0)$ of an $N$-dimensional
manifold and the germ $(L,0)\subset (S,0)$ of an $n$-dimensional
reduced subvariety. $H^\bullet_{Giv}(S,L,0)$
denotes the cohomology of the de Rham complex
\begin{eqnarray}\label{2.9}
(\Omega_{S,0}^\bullet/\{\omega\in\Omega^\bullet_{S,0}\,|\,
\omega|_{L^{reg}}=0\},
\end{eqnarray}
which was considered first by Givental \cite[ch. 1.1]{Gi88}.

\begin{theorem}\label{t2.16}
(a) \cite[ch. 1.1]{Gi88} If $(L,0)$ is quasihomogeneous then 
$H^\bullet_{Giv}(S,L,0)=0$.

\medskip
(b) \cite{Va85} 
If $(S,0)=(\C^2,0)$ and $(L,0)=(f^{-1}(0),0)$ for 
a holomorphic function germ $f:(\C^2,0)\to(\C,0)$
with an isolated singularity at 0, then
$\dim H^1_{Giv}(S,L,0)=\mu-\tau$, where 
\begin{eqnarray*}
\mu:=\dim\OO_{\C^2,0}/\Bigl(\frac{\paa f}{\paa x_i}\Bigr)
\quad\textup{and}\quad
\tau:=\dim\OO_{\C^2,0}/\Bigl(f,\frac{\paa f}{\paa x_i}\Bigr).
\end{eqnarray*}

(c) \cite[ch. 1.2]{Gi88} In the situation of (b),
$\mu>\tau\iff(L,0)$ is not quasihomogeneous.
And if $(L,0)$ is not quasihomogeneous, then 
 $\eta\in\Omega^1_{\C^2,0}$ satisfies 
$[\eta]\in H^1_{Giv}(\C^2,L,0)-\{0\}$ if
$\ddd \eta=u(x_1,x_2)\ddd x_1\ddd x_2$ with $u(0)\neq 0$
(i.e. $\ddd\eta$ is a volume form).

\medskip
(d) \cite[ch 1.1]{Gi88} 
In part (b) in Theorem \ref{t2.15}, $F$ is holomorphic
on $L_M$ if and only if 
$[\alpha]=0\in H^1_{Giv}(T^*M,L_M,(y^0,t^0))$.

\medskip
(e) In part (b) in Theorem \ref{t2.15}, suppose that
$\textup{embdim}(L_M,(y^0,t^0))\leq n+1$. Then 
$(L_M,(y^0,t^0))\cong (\C^{n-1},0)\times (C,0)$ where
$(C,0)$ is the germ of a plane curve. 
And then $F$ is holomorphic on $L$ if and only if $(C,0)$ 
is quasihomogeneous.

\medskip
(f) A germ $(M,t^0)$ of a simple $F$-manifold has a 
(holomorphic) Euler field.
\end{theorem}

{\bf Proof of the parts (e) and (f):} 
(e) The first statement follows from
Theorem \ref{t2.14}, and the decomposition is compatible
with a decomposition of the symplectic germ
$(T^*M,(y^0,t^0))$. The second statement follows from the first
statement and from the parts (a), (c) and (d). 

(f) We can restrict to an irreducible germ $(M,t^0)$ of a
simple $F$-manifold. The caustic $\KK$ is a hypersurface.
At a generic point $p\in\KK$, 
\begin{eqnarray*}
0=\dim(S_\mu(p),p)-l(P(p))=n-1-l(P(p)), \quad\textup{so }
l(P(p))=n-1,
\end{eqnarray*}
so $P(p)=(2,1,...,1)$, and $(M,t^0)$ is a product of $n-2$
1-dimensional and 1 2-dimensional $F$-manifolds.
They have Euler fields, so $F$ is holomorphic on $M-\KK^{sing}$.
But $\codim \KK^{sing}\geq 2$, so $F$ is holomorphic on $M$,
and the Euler field $E_0$ from Theorem \ref{t2.15} extends
to $M$. 
\hfill$\Box$

\bigskip
A generalization of the generically semisimple $F$-manifolds
are the generically regular $F$-manifolds.

\begin{definition}\label{t2.17}\cite[Definition 1.2]{DH17}
Let $(M,\circ,e,E)$ be an $F$-manifold with Euler field.

\medskip
(a) The Euler field is regular at a point $t\in M$
if $E\circ|_t:T_tM\to T_tM$ is a regular endomorphism,
i.e. it has for each eigenvalue only one Jordan block.

\medskip
(b) The $F$-manifold with Euler field $(M,\circ,e,E)$
is called a [generically] regular $F$-manifold 
if the Euler field is regular at all
[respectively at generic] points.
\end{definition}

Theorem 1.3 in \cite{DH17} provides a generalization
of the canonical coordinates of a semisimple
$F$-manifold with Euler field to the case of a 
regular $F$-manifold.

\section{2-dimensional $F$-manifolds}\label{c3}
\setcounter{equation}{0}

The 2-dimensional germs of $F$-manifolds were classified
in \cite{He02}.

\begin{theorem}\label{t3.1}\cite[Theorem 4.7]{He02}
In dimension 2, (up to isomorphism) the germs of
$F$-manifolds fall into three types:

\medskip
(a) The semisimple germ (of type $A_1^2$). 
See Lemma \ref{t2.8} (c) for it 
and for the Euler fields on it.

\medskip
(b) Irreducible germs, which (i.e. some holomorphic 
representatives of them) are at generic points semisimple.
They form a series $I_2(m)$, $m\in\Z_{\geq 3}$.
The germ of type $I_2(m)$ can be given as follows.
\begin{eqnarray}
(M,0)&=&(\C^2,0)\ \textup{with coordinates }t=(t_1,t_2)
\textup{ and }\paa_k:=\frac{\paa}{\paa t_k},\nonumber\\
e&=&\paa_1,\quad \paa_2\circ\paa_2 =t_2^{m-2}e.\label{3.1}
\end{eqnarray}
Any Euler field takes the shape
\begin{eqnarray}\label{3.2}
E&=& (t_1+c_1)\paa_1 + \frac{2}{m}t_2\paa_2
\quad\textup{for some }c_1\in\C.
\end{eqnarray}

\medskip
(c) An irreducible germ, such that the multiplication is
everywhere irreducible. It is called $\NN_2$, and it
can be given as follows.
\begin{eqnarray}
(M,0)&=&(\C^2,0)\ \textup{with coordinates }t=(t_1,t_2)
\textup{ and }
\paa_k:=\frac{\paa}{\paa t_k},\nonumber\\
e&=&\paa_1,\quad \paa_2\circ\paa_2 =0.\label{3.3}
\end{eqnarray}
Any Euler field takes the shape
\begin{eqnarray}\label{3.4}
E&=& (t_1+c_1)\paa_1 + g(t_2)\paa_2
\quad\textup{for some }c_1\in\C\\
&&\hspace*{3cm}\textup{ and some function }
g(t_2)\in\C\{t_2\}.\nonumber
\end{eqnarray}
\end{theorem}

But in the case of $\NN_2$, one has still freedom in the
choice of the coordinate $t_2$, and one can use this to
put an Euler field into a normal form. This was not studied in
\cite{He02}, but in \cite{DH20-1}.

\begin{theorem}\label{t3.2}\cite[Theorem 48]{DH20-1}
(a) The automorphism group of the germ $\NN_2$ of an
$F$-manifold is
\begin{eqnarray}\label{3.5}
&&\Aut(\NN_2)=\Aut((M,0),\circ,e,E)\\
&=&\{(t_1,t_2)\mapsto (t_1,f(t_2))\,|\, f(t_2)\in\C\{t_2\}
\textup{ with }f(0)=0,f'(0)\neq 0\}.\nonumber
\end{eqnarray}

(b) Let $\www E$ be an Euler field on $\NN_2$. 
Its orbit under the automorphism group $\Aut(\NN_2)$ contains
precisely one of the Euler fields in the following list,
\begin{eqnarray}\label{3.6}
&&E=(t_1+c)\paa_1+\paa_2,\\
&&E=(t_1+c)\paa_1,\label{3.7}\\
&&E=(t_1+c)\paa_1+c_0t_2\paa_2,\label{3.8}\\
&&E=(t_1+c)\paa_1+t_2^r(1+c_1t_2^{r-1})\paa_2,\label{3.9},
\end{eqnarray}
where $c,c_1\in\C$, $c_0\in\C^*$ and $r\in\Z_{\geq 2}$. 
\end{theorem}

$\NN_2$ with the Euler field $E$ in \eqref{3.6} is regular
(Definition \ref{t2.17}). $\NN_2$ with the Euler field
$E$ in \eqref{3.8} or \eqref{3.9} is generically regular.
$\NN_2$ with the Euler field in \eqref{3.7} is not even
generically regular.

\section{Basic formulas for 3-dimensional $F$-manifolds}
\label{c4}
\setcounter{equation}{0}

\begin{notations}\label{t4.1}
In the sections \ref{c4} and \ref{c5}, we consider a 
3-dimensional complex manifold $M$ with a holomorphic 
commutative multiplication on the holomorphic
tangent bundle and with a unit field $e$ 
(so $e\circ=\id$) with $\textup{Lie}_e(\circ)=0$. 
This condition $\textup{Lie}_e(\circ)=0$ follows from
the integrability condition \eqref{1.1} of an $F$-manifold, 
which we do not suppose at the beginning. 
But it would be a nuisance not to suppose 
$\textup{Lie}_e(\circ)=0$ from the beginning.

We work locally near a point $p\in M$ 
and suppose to have coordinates 
$t=(t_1,t_2,t_3)$ with $\paa_1=e$ and 
$(M,p)\cong(\C^3,0)$ and 
coordinate vector fields $\paa_j=\paa/\paa t_j$. 
Let $y=(y_1,y_2,y_3)$ be the fiber coordinates on $T^*M$ 
which correspond to $\paa_1,\paa_2,\paa_3$. 
Then the canonical 1-form $\alpha$ takes the shape 
$\alpha=\sum_{i=1}^3y_i\ddd t_i$.

We write 
\begin{eqnarray}\label{4.1}
\paa_2\circ\paa_2 &=& \www a_1\paa_1+\www a_2\paa_2+a_3\paa_3,\\
\paa_2\circ\paa_3 &=& \www b_1\paa_1+b_2\paa_2+b_3\paa_3,
\label{4.2}\\
\paa_3\circ\paa_3 &=& \www c_1\paa_1+c_2\paa_2+\www c_3\paa_3,
\label{4.3}
\end{eqnarray}
with $\www{a}_1,\www{a}_2,a_3,\www{b}_1,b_2,b_3,
\www{c}_1,c_2,\www{c}_3\in \OO_M$. 
Many formulas take a simpler shape 
if we rewrite the formulas above as follows,
\begin{eqnarray}\label{4.4}
(\paa_2-b_3\paa_1)\circ(\paa_2-b_3\paa_1) 
&=& a_1\paa_1+a_2(\paa_2-b_3\paa_1)+a_3(\paa_3-b_2\paa_1),\\
(\paa_2-b_3\paa_1)\circ(\paa_3-b_2\paa_1) 
&=& b_1\paa_1, \label{4.5}\\
(\paa_3-b_2\paa_1)\circ(\paa_3-b_2\paa_1) 
&=& c_1\paa_1+c_2(\paa_2-b_3\paa_1)+ c_3(\paa_3-b_2\paa_1),
\hspace*{1cm}
\label{4.6}
\end{eqnarray}
with $a_j,b_j,c_j\in\OO_M$. The condition 
$\textup{Lie}_e(\circ)=0$ is equivalent to 
$a_j,b_j,c_j\in\C\{t_2,t_3\}$, so we suppose this from
now on (and it also implies 
$\www a_1,\www a_2,\www b_1,\www c_1,\www c_3\in\C\{t_2,t_3\}$).
We denote 
\begin{eqnarray}\label{4.7}
\paa_ia_j:=a_{ji},\ \paa_ib_j:=b_{ji},\ \paa_ic_j:=c_{ji},
\textup{ and analogously for }\www a_j,\www b_j,\www c_j.
\end{eqnarray}

If $s=(s_1,s_2,s_3)$ is another system of coordinates on
$(M,p)$ with $t=t(s)$ and $s=s(t)$ and $s(p)=0$, 
write $\www\paa_j:=\paa/\paa s_j$
for the coordinate vector fields of this system of 
coordinates, and write $z=(z_1,z_2,z_3)$ for the 
fiber coordinates which correspond to
$\www\paa_1,\www\paa_2,\www\paa_3$. Then
\begin{eqnarray}\label{4.8}
\ddd t_i= \sum_{j=1}^3\www\paa_jt_i\cdot \ddd s_j,\quad 
z_j=\sum_{i=1}^3\www\paa_j t_i\cdot y_i.
\end{eqnarray} 
We suppose $\www\paa_1=e=\paa_1$. This is equivalent to 
$t_i(s)\in(\delta_{i1}\cdot s_1+\C\{s_2,s_3\})$
and also to 
$s_j(t)\in(\delta_{j1}\cdot t_1+\C\{t_2,t_3\})$.
Often it is useful to make first a special 
coordinate change of the type 
$t_2=s_2,t_3=s_3,t_1=s_1+\tau$ with
$\tau\in\C\{t_2,t_3\}=\C\{s_2,s_3\}$. Then
\begin{eqnarray}\label{4.9}
z_1=y_1,\quad z_2=\paa_2\tau\cdot y_1+y_2,\quad
z_3=\paa_3\tau\cdot y_1+y_3.
\end{eqnarray}
\end{notations}

\begin{lemma}\label{t4.2}
In the situation of the Notations \ref{t4.1}, 
the multiplication is associative if and only if
\begin{eqnarray}\label{4.10}
a_1=-a_3c_3,\quad b_1=a_3c_2,\quad c_1=-a_2c_2.
\end{eqnarray}
\end{lemma}

{\bf Proof:} Straightforward calculations with 
\eqref{4.4}--\eqref{4.6} of both sides of the
equations
\begin{eqnarray*}
&&\bigl((\paa_2-b_3\paa_1)\circ(\paa_2-b_3\paa_1)\bigr)\circ(\paa_3-b_2\paa_1)\\
&=&\bigl((\paa_2-b_3\paa_1)\circ(\paa_3-b_2\paa_1)\bigr)\circ(\paa_2-b_3\paa_1),\\
&&\bigl((\paa_2-b_3\paa_1)\circ(\paa_3-b_2\paa_1)\bigr)\circ(\paa_3-b_2\paa_1)\\
&=&\bigl((\paa_3-b_2\paa_1)\circ(\paa_3-b_2\paa_1)\bigr)\circ(\paa_2-b_3\paa_1).
\end{eqnarray*}
\hfill$\Box$

The next lemma starts with $M$ as in the Notations \ref{t4.1},
but with associative multiplication, 
and tells to which of the four algebras
$Q^{(1)},Q^{(2)},Q^{(3)}$ or $Q{(4)}$ in the Remarks 
\ref{t1.1} the algebra $T_tM$ for $t\in M$ is isomorphic.

\begin{lemma}\label{t4.3}
Let $(M,\circ,e)$ be as in the Notations \ref{t4.1}
with coordinates $t=(t_1,t_2,t_3)$, 
and suppose that the multiplication $\circ$ is associative.
Define $R_1,R_2,R_3\in\C\{t_2,t_3\}$ by 
\begin{eqnarray}\label{4.11}
R_1:= a_3c_3-\frac{1}{3}a_2^2,\
R_2:= a_2c_2-\frac{1}{3}c_3^2,\
R_3:= a_3c_2-\frac{1}{9}a_2c_3.
\end{eqnarray}
For a point $t\in M$ the following statements hold.
\begin{eqnarray}\label{4.12}
T_tM\cong Q^{(1)}&\iff& (a_2,a_3,c_2,c_3)(t)=0.\\
T_tM\cong Q^{(2)}&\iff& (R_1,R_2,R_3)(t)=0,\ 
(a_3,c_2)(t)\neq 0.\label{4.13}\\
T_tM\cong Q^{(3)}&\iff& (9R_3^2-4R_1R_2)(t)=0, 
(R_1,R_2,R_3)(t)\neq 0.\hspace*{1cm}\label{4.14}\\
T_tM\cong Q^{(4)}&\iff& (9R_3^2-4R_1R_2)(t)\neq 0.
\label{4.15}\\
a_3(t)\neq 0&\textup{and}&(R_1,R_3)(t)=0
\ \Rightarrow \ R_2(t)=0.
\label{4.16}\\
c_2(t)\neq 0&\textup{and}&(R_2,R_3)(t)=0
\ \Rightarrow\ R_1(t)=0.
\label{4.17}
\end{eqnarray}
\end{lemma}

{\bf Proof:}
Define 
\begin{eqnarray}\label{4.18}
\psi_1:=\paa_2-b_3\paa_1-\frac{1}{3}a_2\paa_1\ \textup{ and }\ 
\psi_2:=\paa_3-b_2\paa_1-\frac{1}{3}c_3\paa_1.
\end{eqnarray}
One calculates
\begin{eqnarray}\label{4.19}
\psi_1^{\circ 2}&=&\frac{1}{3}a_2\psi_1 +a_3\psi_2
-\frac{2}{3}R_1\paa_1,\\
\psi_1\circ\psi_2&=& -\frac{1}{3}c_3\psi_1-\frac{1}{3}a_2\psi_2
+R_3\paa_1,\label{4.20}\\
\psi_2^{\circ 2}&=&c_2\psi_1 +\frac{1}{3}c_3\psi_2
-\frac{2}{3}R_2\paa_1,\label{4.21}
\end{eqnarray}
\begin{eqnarray}\label{4.22}
0&=&\psi_1^{\circ 3}+R_1\psi_1+
(\frac{2}{9}a_2R_1-a_3R_3)\paa_1,\\
0&=&\psi_2^{\circ 3}+R_2\psi_2+
(\frac{2}{9}c_3R_2-c_2R_3)\paa_1.\label{4.23}
\end{eqnarray}
The lack of a quadratic term $\psi_1^{\circ 2}$ in \eqref{4.22}
shows that the sum of the three eigenvalues of 
$\psi_1\circ|_t:T_tM\to T_tM$ is zero for any $t\in M$, 
and similarly for
$\psi_2\circ|_t$. If $T_tM$ is irreducible, then 
$\psi_1\circ|_t$ and $\psi_2\circ|_t$ have only one 
eigenvalue, which is then  zero. 
Therefore they are nilpotent, so $\psi_1^{\circ 3}|_t=0$
and $\psi_2^{\circ 3}|_t=0$, so $(R_1,R_2,R_3)(t)=0$.
Vice versa, if $(R_1,R_2,R_3)(t)=0$, then
$\psi_1^{\circ 3}|_t=0$ and $\psi_2^{\circ 3}|_t=0$,
so $\psi_1\circ|_t$ and $\psi_2\circ|_t$
are nilpotent, and $T_tM$ is an irreducible algebra.
We proved
\begin{eqnarray*}
T_tM\cong Q^{(1)}\textup{ or }T_tM\cong Q^{(2)}
\iff (R_1,R_2,R_3)(t)=0.
\end{eqnarray*}
Suppose that $T_tM$ is an irreducible algebra.
If $a_3(t)\neq 0$ then 
\begin{eqnarray*}
(\paa_3-b_2\paa_1)(t)=a_3(t)^{-1}
\Bigl((\paa_2-b_3\paa_1)^{\circ 2}
-a_2(\paa_2-b_3\paa_1)-a_1\paa_1\Bigr)(t),\\
\textup{so }\quad 
T_tM=\bigoplus_{j=0}^2\C\cdot \paa_2^{\circ j}(t),
\quad\textup{ and thus }\quad T_tM\cong Q^{(2)},
\end{eqnarray*}
and in the same way $c_2(t)\neq 0$ implies 
$T_tM\cong Q^{(2)}$. The other way round,
if $(a_3,c_2)(t)=0$, then $(R_1,R_2)(t)=0$ implies
also $(a_2,c_3)(t)=0$, and then $T_tM\cong Q^{(1)}$. 
This finishes the proof of \eqref{4.12} and \eqref{4.13}.

Next we want to show \eqref{4.14}. 
\eqref{4.22} and \eqref{4.23} generalize as follows
for arbitrary $\psi:=\lambda_1\psi_1+\lambda_2\psi_2$
with $\lambda_1,\lambda_2\in\C\{t_2,t_3\}$:
\begin{eqnarray}\label{4.24}
0&=&\psi^{\circ 3} +
[R_1\lambda_1^2-3R_3\lambda_1\lambda_2+R_2\lambda_2^2]
\cdot\psi\\
&+&\Bigl[(\frac{2}{9}a_2R_1-a_3R_3)\lambda_1^3
-(\frac{2}{3}c_3R_1-a_2R_3)\lambda_1^2\lambda_2\nonumber\\
&&-(\frac{2}{3}a_2R_2-c_3R_3)\lambda_1\lambda_2^2
+(\frac{2}{9}c_3R_2-c_2R_3)\lambda_2^3\Bigr]\cdot\paa_1.
\nonumber
\end{eqnarray}
A lengthy calculation shows that the 
discriminant of \eqref{4.24} is 
\begin{eqnarray}\label{4.25}
&&4(\textup{coefficient of }\psi)^3+27(\textup{coefficient of }
\paa_1)^2\\
&=&(9R_3^2-4R_1R_2)\cdot 3(a_3\lambda_1^3-
a_2\lambda_1^2\lambda_2+c_3\lambda_1\lambda_2^2
-c_2\lambda_2^3)^2.\nonumber
\end{eqnarray}
First suppose $T_tM\cong Q^{(3)}$. Then $\psi|_t$ has at most 
two different eigenvalues, and the discriminant in 
\eqref{4.25} must vanish at $t$.  Because $\lambda_1$
and $\lambda_2\in\C\{t_2,t_3\}$ are arbitrary and 
$(a_2,a_3,c_2,c_3)\neq 0$, this shows 
$(9R_3^2-4R_1R_2)(t)=0$. Vice versa, if 
$(9R_3^2-4R_1R_2)(t)=0$, then the discriminant in \eqref{4.25}
vanishes at $t$ for any $\psi$. Therefore $\psi\circ|_t$
has at most two eigenvalues for any $\psi$.
This shows $T_tM\not\cong Q^{(4)}$. 
Then the condition $(R_1,R_2,R_3)\neq 0$ yields
$T_tM\cong Q^{(3)}$. 
This proves \eqref{4.14}. 

\eqref{4.15} is a 
consequence of \eqref{4.12}--\eqref{4.14}.
The implications \eqref{4.16} and \eqref{4.17}
are trivial.\hfill$\Box$

\begin{lemma}\label{t4.4}
Let $(M,\circ,e)$ be as in the Notations \ref{t4.1}
with coordinates $t=(t_1,t_2,t_3)$ and fiber coordinates
$y=(y_1,y_2,y_3)$ of $T^*M$, and suppose that the multiplication
$\circ$ is associative. With the result from Lemma \ref{t4.2}, 
the ideal $\II_M\subset \OO(T^*M)$ which defines 
the analytic spectrum is
\begin{eqnarray}\label{4.26}
\II_M&=&\bigl(y_1-1,Y_{22},Y_{23},Y_{33}\bigr),\quad
\textup{where}\\
Y_{22}&:=&(y_2-b_3)(y_2-b_3) 
+a_3c_3-a_2(y_2-b_3)-a_3(y_3-b_2),\nonumber\\
Y_{23}&:=&(y_2-b_3)(y_3-b_2) -a_3c_2,\nonumber\\
Y_{33}&:=&(y_3-b_2)(y_3-b_2)
+a_2c_2-c_2(y_2-b_3)-c_3(y_3-b_2),\nonumber
\end{eqnarray}
with $a_j,b_j,c_j\in\OO_M$. 
Recall the notation in formula \eqref{4.7}. Define
\begin{eqnarray}
A_2&:=& a_2(-b_{22}+b_{33}+a_{23})
+a_3(-2c_{22}-c_{33})-a_{32}c_2-a_{33}c_3,\nonumber\\
A_2^{dual}&:=&c_3(-b_{33}+b_{22}+c_{32})
+c_2(-2a_{33}-a_{22})-c_{23}a_3-c_{22}a_2,\nonumber\\
A_3&:=& -3b_{22}+3b_{33}+a_{23}-c_{32}.\label{4.27}
\end{eqnarray}
Then
\begin{eqnarray}
\{y_1-1,Y_{ij}\} &=& \paa_1 Y_{ij} \quad (=0\ \textup{ because of }
a_k,b_k,c_k\in\C\{t_2,t_3\}),\nonumber\\ 
\{Y_{22},Y_{23}\}&=& Y_{22} [-2b_{22}+2b_{33}+a_{23}]
+ Y_{23} [a_{22}+a_{33}] 
+ Y_{33} [a_{33}]\nonumber\\
&+& (y_2-b_3) A_2 
+ (y_3-b_2) a_3A_3 
+ [-a_3A_2^{dual}-a_3c_3A_3],\nonumber\\
\{Y_{33},Y_{23}\}&=& Y_{33} [-2b_{33}+2b_{22}+c_{32}]
+ Y_{23} [c_{33}+c_{22}] 
+ Y_{22} [c_{22}]\nonumber\\
&+& (y_3-b_2) A_2^{dual} 
- (y_2-b_3) c_2A_3 
+[-c_2A_2+c_2a_2A_3],\nonumber\\
\{Y_{22},Y_{33}\} &=& 
Y_{22} [-2c_{22}] 
+ Y_{23} [2(-2b_{22}+2b_{33}+a_{23}-c_{32})] \nonumber\\
&+& Y_{33} [2a_{33}]
+ (y_2-b_3) [-A_2^{dual}-c_3A_3] \nonumber\\
&+& (y_3-b_2) [A_2-a_2A_3] \nonumber\\
&+& [-c_3A_2+a_2A_2^{dual}+(a_2c_3+a_3c_2)A_3].\label{4.28}
\end{eqnarray}
Therefore, $(M,\circ,e)$ is an $F$-manifold if and only if
\begin{eqnarray}\label{4.29}
(a_2,a_3,c_2,c_3)=0,\\
\textup{or}\qquad (A_2,A_2^{dual},A_3)=0.\label{4.30}
\end{eqnarray}
The intersection of these two cases is the case 
\eqref{4.29} with additionally $b_{22}-b_{33}=0$.
\end{lemma}

{\bf Proof:}
The calculation of the Poisson brackets in \eqref{4.28}
is straightforward and leads to the claimed formulas in
\eqref{4.28}. 
By Theorem \ref{t2.12}, $(M,\circ,e)$ is an $F$-manifold 
if and only if $\{\II_M,\II_M\}\subset \II_M$, so
if and only if 
\begin{eqnarray}\label{4.31}
A_2=A_2^{dual}=a_3A_3=c_2A_3=c_3A_3=a_2A_3=0.
\end{eqnarray}
This leads to the two cases \eqref{4.29} and \eqref{4.30}.
\hfill$\Box$

\begin{remark}\label{t4.5}
In the case \eqref{4.30}, the condition $A_3=0$ can be used
to make a specific special coordinate change
as in \eqref{4.9}, namely we choose the new coordinates
$s=(s_1,s_2,s_3)$ such that
\begin{eqnarray}\label{4.32}
t_2=s_2,t_3=s_3,t_1=s_1+\tau\quad\textup{with}\quad
\tau\in\C\{t_2,t_3\}\textup{ with }\\
\paa_2\tau=-b_3-\frac{1}{3}a_2,
\ \paa_3\tau=-b_2-\frac{1}{3}c_3.\nonumber
\end{eqnarray}
With the notation $\www\paa_j:=\paa/\paa s_j$ 
and with $\psi_1,\psi_2$ as in \eqref{4.18}, 
we obtain
\begin{eqnarray}\label{4.33}
\www\paa_1=\paa_1=e,\ 
\www\paa_2=\psi_1,\ \www\paa_3=\psi_2,\\
\paa_2-b_3\paa_1=\www\paa_2+\frac{1}{3}a_2\paa_1,\quad  
\paa_3-b_2\paa_1=\www\paa_3+\frac{1}{3}c_3\paa_1.\nonumber 
\end{eqnarray}
Formula \eqref{4.24} in the proof of Lemma \ref{t4.4} 
tells that for any 
$\psi=\lambda_1\psi_1+\lambda_2\psi_2$ with 
$\lambda_1,\lambda_2\in\C\{t_2,t_3\}$, the sum of the
eigenvalues of $\psi\circ$ is zero at any $t\in M$. 

If we call now the new coordinates again $t=(t_1,t_2,t_3)$,
the new and old coefficients $a_2,a_3,c_2,c_3$ coincide,
and the new coefficients $b_2^{(new)},b_3^{(new)},A_2^{(new)},
(A_2^{dual})^{(new)},A_3^{(new)}$ become
\begin{eqnarray}\label{4.34}
b_3^{(new)}=-\frac{1}{3}a_2,\ b_2^{(new)}=-\frac{1}{3}c_3,\ 
A_3^{(new)}=0,\\
A_2^{(new)}=-\paa_3R_1 + \frac{1}{3}a_2c_{32}-2a_3c_{22}
-a_{32}c_2,\label{4.35}\\
(A_2^{dual})^{(new)}=-\paa_2R_2 + 
\frac{1}{3}a_{23}c_3-2a_{33}c_2 -a_3c_{23}.\label{4.36}
\end{eqnarray}
In the case \eqref{4.30}, we will often, but not always,
assume that the coordinates $t=(t_1,t_2,t_3)$ have
been chosen as in this remark.
\end{remark}

If a germ $(M,0)$ of a 3-dimensional 
$F$-manifold satisfies $T_0M\cong Q^{(2)}$ or $Q^{(3)}$
or $Q^{(4)}$, 
then life is easier than in the case $T_0M\cong Q^{(1)}$.
The next lemma makes this explicit in one way.

\begin{lemma}\label{t4.6}
Let $((M,0),\circ,e)$ be as in the Notations \ref{t4.1}
with coordinates $t=(t_1,t_2,t_3)$ with $t(0)=0$ 
and fiber coordinates $(y_1,y_2,y_3)$ of $T^*M$, 
and suppose that the multiplication
$\circ$ is associative. Suppose $T_0M\not\cong Q^{(1)}$. 
The coordinates $t$ can and will be chosen such that
\begin{eqnarray}\label{4.37}
\C\cdot\paa_1|_0\oplus\C\cdot \paa_2|_0\oplus \C\cdot 
\paa_2^{\circ 2}|_0=T_0M.
\end{eqnarray}
Then
\begin{eqnarray}\label{4.38}
\paa_2^{\circ 3}&=& g_2\cdot\paa_2^{\circ 2}+g_1\cdot\paa_2
+g_0\cdot\paa_1,\\
\paa_3 &=& h_2\cdot \paa_2^{\circ 2}+h_1\cdot \paa_2
+h_0\cdot\paa_1\label{4.39}
\end{eqnarray}
for suitable coefficients $g_2,g_1,g_0,h_2,h_1,h_0
\in\C\{t_2,t_3\}$. 
We denote similarly to \eqref{4.7}
\begin{eqnarray*}
\paa_ig_j=:g_{ji},\quad \paa_ih_j=:h_{ji}.
\end{eqnarray*}
The ideal $\II_M\subset\OO_{T^*M}$ which defines the analytic
spectrum is
\begin{eqnarray}\label{4.40}
\II_M&=& (y_1-1,Z_2,Z_3\},\quad\text{where}\\
Z_2&:=& y_2^3-g_2y_2^2-g_1y_2-g_0,\nonumber\\
Z_3&:=& y_3-h_2y_2^2-h_1y_2-h_0.\nonumber
\end{eqnarray}
Then 
\begin{eqnarray}
&&\{y_1-1,Z_j\}=\paa_1Z_j\quad(=0\textup{ because of }
g_i,h_i\in\C\{t_2,t_3\}),\nonumber\\
&&\{Z_3,Z_2\}= Z_2[2g_{22}h_2+(3y_2+g_2)h_{22}
+3h_{12}]\nonumber\\
&+& y_2^2[\paa_2((g_2^2+2g_1)h_2+g_2h_1+3h_0)-g_{23}]
\label{4.41}\\
&+& y_2[(2g_{22}g_1+2g_{02})h_2+(g_2g_1+3g_0)h_{22}+g_{12}h_1
\nonumber\\
&& \hspace*{4cm} +2g_1h_{12}-2g_2h_{02}-g_{13}]\label{4.42}\\
&+& [2g_{22}g_0h_2+g_2g_0h_{22}+g_{02}h_1+3g_0h_{12}
-g_1h_{02}-g_{03}].\label{4.43}
\end{eqnarray}
Therefore $((M,0),\circ,e)$ is a germ of an $F$-manifold
if and only if the terms in square brackets in
\eqref{4.41}--\eqref{4.43} vanish. 
\end{lemma}

{\bf Proof:}
In each of the algebras $Q^{(j)}$ for $j\in\{2,3,4\}$, 
a generic element $a$ satisfies $Q^{(j)}=\C\cdot 1\oplus
\C\cdot a\oplus \C\cdot a^{\circ 2}$. One can choose
the coordinates $t$ on $(M,0)$ such that $\paa_1=e$
and $\paa_2|_0$ is such a generic element. 
This implies \eqref{4.37}--\eqref{4.40}. 
The calculation of $\{Z_3,Z_2\}$ is straightforward.
\hfill$\Box$

\begin{corollary}\label{t4.7}
Let $g_2^{(0)},g_1^{(0)},g_0^{(0)}\in\C\{t_2\}$
and $h_2,h_1,h_0\in\C\{t_2,t_3\}$ be arbitrary.
There exist unique $g_2,g_1,g_0\in\C\{t_2,t_3\}$
such that $g_j|_{t_3=0}=g_j^{(0)}$ and such that
the 3-dimensional germ $(M,0)$ of a manifold with
multiplication $\circ$ on $TM$ defined by $\paa_1=e$, 
\eqref{4.38} and \eqref{4.39} is a germ of an $F$-manifold.
\end{corollary}

{\bf Proof:} The Cauchy-Kovalevski theorem in the following
form \cite[(1.31), (1.40), (1.41)]{Fo95} will be applied
(there the setting is real analytic, but the proofs and 
statements hold also in the complex analytic setting):
Given $N\in\N$ and matrices 
$A_i,B\in M_{N\times N}(\C\{s_1,...,s_m,y,x_1,...,x_N\})$,
there exists a unique vector 
\begin{eqnarray*}
\Phi\in M_{N\times 1}(\C\{s_1,...,s_m,y\})
\end{eqnarray*}
with
\begin{eqnarray}\label{4.44}
\frac{\paa\Phi}{\paa y}&=&\sum_{i=1}^mA_i(s,y,\Phi)
\frac{\paa\Phi}{\paa s_i}+B(s,y,\Phi),\\
\Phi(s,0)&=&0.\nonumber
\end{eqnarray}
In our situation $y=t_3$, $(s_1,...,s_m)=(t_2)$, 
$\Phi=(g_2-g_2^{(0)},g_1-g_1^{(0)},g_0-g_0^{(0)})^t$, 
and $A_1,A_2,A_3$ and $B$ come from the terms in 
\eqref{4.41}--\eqref{4.43} without $g_{23},g_{13},g_{03}$,
more precisely, \eqref{4.44} is here
\begin{eqnarray}\label{4.45}
\paa_3\begin{pmatrix}g_2\\g_1\\g_0\end{pmatrix}
= \begin{pmatrix}\paa_2((g_2^2+2g_1)h_2+g_2h_1+3h_0)\\
(2g_{22}g_1+2g_{02})h_2+(g_2g_1+3g_0)h_{22}+g_{12}h_1
\\ \hspace*{2cm} +2g_1h_{12}-2g_2h_{02}\\
2g_{22}g_0h_2+g_2g_0h_{22}+g_{02}h_1+3g_0h_{12}-g_1h_{02}
\end{pmatrix}.
\end{eqnarray}
The Cauchy-Kovalevski theorem tells that there exist unique
$g_2,g_1,g_0\in\C\{t_2,t_3\}$ such that $g_j|_{t_3=0}=g_j^{(0)}$
and such that the terms in 
\eqref{4.41}--\eqref{4.43} vanish. 
The multiplication on $TM$ which is defined by \eqref{4.38}
and \eqref{4.39}, is automatically associative.
The condition $\{\II_M,\II_M\}\subset \II_M$ is equivalent
to the vanishing of the terms in \eqref{4.41}--\eqref{4.43}.
By Theorem \ref{t2.12}, $((M,0),\circ,e)$ is an $F$-manifold
if and only if $\{\II_M,\II_M\}\subset \II_M$.
\hfill$\Box$

\begin{remarks}\label{t4.8}
(i) The corollary makes it easy to construct a 3-dimensional
germ $(M,0)$ of an $F$-manifold. Arbitrary initial data
$g_2^{(0)},g_1^{(0)},g_0^{(0)}\in\C\{t_2\}$
and $h_2,h_1,h_0\in\C\{t_2,t_3\}$ give a unique germ
of an $F$-manifold. But this approach does not tell
easily which properties such constructed $F$-manifolds have, 
which families exist and what are their parameters. 

\medskip
(ii) The condition $A_3=0$ in Lemma \ref{t4.4} 
corresponds to the first line of equation \eqref{4.45}.
And we can make a coordinate change as in Remark \ref{t4.5}. 
In the case of an $F$-manifold in Lemma \ref{t4.6},
we can choose new coordinates $s=(s_1,s_2,s_3)$ such that
\begin{eqnarray}\label{4.46}
t_2=s_2,t_3=s_3,t_1=s_1+\tau\quad\textup{with}\quad
\tau\in\C\{t_2,t_3\}\textup{ with}\\
\paa_2\tau=-\frac{1}{3}g_2,
\ \paa_3\tau=-\frac{1}{3}((g_2^2+2g_1)h_2+g_2h_1+3h_0),
\nonumber\\
\www\paa_2=\paa_2+\paa_2\tau\cdot\paa_1,
\ \www\paa_3=\paa_3+\paa_3\tau\cdot\paa_1.
\nonumber
\end{eqnarray}
If we now call the new coordinates again $t=(t_1,t_2,t_3)$, the
new coefficients $g_j^{(new)}$ and $h_j^{(new)}$ satisfy
$g_2^{(new)}=0$ and $2g_1^{(new)}h_2^{(new)}+3h_0^{(new)}=0$.
This says that the sum of the eigenvalues of $\paa_2\circ$
is zero, and that the sum of the eigenvalues of 
$\paa_3\circ$ is 0. The last statement follows from the facts
that the sum of the eigenvalues of $h_0\paa_1\circ$ is
$3h_0$ and that the sum of the eigenvalues of 
$h_2\paa_2^{\circ 2}\circ$ is $2h_2g_1$, because
$\lambda_1^2+\lambda_2^2+\lambda_3^2=-2(\lambda_1\lambda_2
+\lambda_1\lambda_3+\lambda_2\lambda_3)$ for any 
$\lambda_1,\lambda_2,\lambda_3\in\C$ with 
$\sum_{i=1}^3\lambda_i=0$. 
\end{remarks}

\section{3-dimensional not generically semisimple $F$-manifolds}
\label{c5}
\setcounter{equation}{0}

\begin{remarks}\label{t5.1}
We want to classify all 3-dimensional germs of 
$F$-manifolds. The reducible ones are products of
1- and 2-dimensional germs of $F$-manifolds by
Theorem \ref{t2.5}. Those are classified in Lemma \ref{t2.6}
and Theorem \ref{t3.1}. The 3-dimensional reducible
germs of $F$-manifolds are $A_1^3$,
$A_1I_2(m)$ for $m\geq 3$ and $A_1\NN_2$,
and the Euler fields are as described in Theorem \ref{t2.5},
Lemma \ref{t2.6} and Theorem \ref{t3.2}.
It remains to classify the irreducible germs $(M,0)$ of 
$F$-manifolds, i.e. those where $T_0M$ is irreducible. 
We start with those with $T_tM\cong Q^{(1)}$ for any 
$t\in M$. 
\end{remarks}

\begin{theorem}\label{t5.2}
Any 3-dimensional germ $(M,0)$ of an $F$-manifold
with $T_tM\cong Q^{(1)}$ for all $t\in M$ can be given as 
follows:
\begin{eqnarray}
(M,0)&=& (\C^3,0)\textup{ with coordinates }t=(t_1,t_2,t_3),
\nonumber\\
e&=&\paa_1,\ \paa_2^{\circ 2}=
\paa_2\circ(\paa_3-b_2\paa_1)=
(\paa_3-b_2\paa_1)^{\circ 2}=0,\nonumber\\
&&\textup{where }b_2\textup{ is arbitrary in }t_2\C\{t_2,t_3\}.
\label{5.1}
\end{eqnarray}
The ideal $(\paa_2b_2)\subset\C\{t_2,t_3\}$ is up to coordinate
changes an invariant of the germ of an $F$-manifold.
A vector field $E=\varepsilon_1\paa_1+\varepsilon_2\paa_2
+\varepsilon_2\paa_3$ with $\varepsilon_1,\varepsilon_2,
\varepsilon_3\in\C\{t_1,t_2,t_3\}$ is an Euler field if 
and only if $\varepsilon_1\in t_1+\C\{t_2,t_3\}$, 
$\varepsilon_2,\varepsilon_3\in\C\{t_2,t_3\}$ and
\begin{eqnarray}\label{5.2}
\paa_2(\varepsilon_1)=-b_2\paa_2(\varepsilon_2),\quad
\paa_3(\varepsilon_1)=-\varepsilon_2\paa_2(b_2)
-\paa_3(\varepsilon_3b_2)+b_2.
\end{eqnarray}
\end{theorem}

{\bf Proof:} 
Let $((M,0),\circ,e)$ be a germ of an $F$-manifold with
$T_tM\cong Q^{(1)}$ for any $t\in M$. 
We choose coordinates $t=(t_1,t_2,t_3)$ with $t(0)=0$
and use the Notations \ref{t4.1} and Lemma \ref{t4.2}.
By Lemma \ref{t4.3}, $(a_3,a_2,c_3,c_2)=0$. 
We are in the case \eqref{4.29} in Lemma \ref{t4.4}.
The $F$-manifold condition gives no constraint 
on $b_2,b_3\in\C\{t_2,t_3\}$. 

We make a specific special coordinate change
as in \eqref{4.9}, namely we choose the new coordinates
$s=(s_1,s_2,s_3)$ such that
\begin{eqnarray}\label{5.3}
t_2=s_2,t_3=s_3,t_1=s_1+\tau\quad\textup{with}\quad
\tau\in\C\{t_2,t_3\}\textup{ with }\\
\paa_2\tau=-b_3,
\ \paa_3\tau+b_2\in t_2\C\{t_2,t_3\}.\nonumber
\end{eqnarray}
Then $\tau$ exists and is unique. 
With the notation $\www\paa_j:=\paa/\paa s_j$ 
we obtain
\begin{eqnarray}\label{5.4}
\www\paa_1=\paa_1=e,\ 
\www\paa_2=\paa_2-b_3\paa_1,\ 
\www\paa_3=\paa_3+\paa_3\tau\cdot \paa_1,\\
\paa_3-b_2\paa_1=\www\paa_3-(\paa_3\tau+b_2\paa_1).\nonumber 
\end{eqnarray}
If we call now the new coordinates again $t=(t_1,t_2,t_3)$,
the new coefficients $b_2^{(new)}$ and $b_3^{(new)}$ are
\begin{eqnarray*}
b_3^{(new)}=0\quad\textup{and}\quad 
b_2^{(new)}=\paa_3\tau+b_2\in t_2\C\{t_2,t_3\}.
\end{eqnarray*}

Now we want to show that the ideal 
$(\paa_2b_2)\subset\C\{t_2,t_3\}$
is up to coordinate changes an invariant of the germ of
an $F$-manifold. 
We consider a coordinate change as in \eqref{4.8}
with $z_1=y_1$, which implies $\www\paa_1t_1=1$ and 
$\www\paa_1t_2=\www\paa_1t_3=0$. Then 
\begin{eqnarray*}
z_2&=&\www\paa_2t_1\cdot y_1+\www\paa_2t_2\cdot y_2
+\www\paa_2t_3\cdot y_3,\\
z_3&=&\www\paa_3t_1\cdot y_1+\www\paa_3t_2\cdot y_2
+\www\paa_3t_3\cdot y_3.
\end{eqnarray*}
$z_2-(\www\paa_2t_1+\www\paa_2t_3\cdot b_2)y_1$ and 
$z_3-(\www\paa_3t_1+\www\paa_3t_3\cdot b_2)y_1$ are
nilpotent in $\OO_{T^*M}/\II_M|_{L_M}$. 
We need $z_2$ to be nilpotent, so the coordinate change
satisfies
\begin{eqnarray*}
\www\paa_2t_1=-\www\paa_2t_3\cdot b_2.
\end{eqnarray*}
And then 
\begin{eqnarray*}
\www b_2 := \www\paa_3t_1+\www\paa_3t_3\cdot b_2.
\end{eqnarray*}
takes the role of $b_2$ for the new coordinates. 
A short calculation shows
\begin{eqnarray*}
\www\paa_2\www b_2= (\paa_2 b_2)(t(s))\cdot
(\www\paa_3t_3\cdot\www\paa_2t_2-\www\paa_2t_3\cdot
\www\paa_3t_2).
\end{eqnarray*}
The second factor is a unit.
Therefore the ideal $(\paa_2b_2)$ is up to coordinate
changes an invariant of the germ $(M,0)$ of the
$F$-manifold.

The constraint \eqref{5.2} for an Euler field
$E=\varepsilon_1\paa_1+\varepsilon_2\paa_2+\varepsilon_3\paa_3$
follows straightforwardly from the explicit version
\begin{eqnarray}\label{5.5}
0=[E,\paa_i\circ\paa_j]-[E,\paa_i]\circ \paa_j - 
[E,\paa_j]\circ\paa_i - \paa_i\circ\paa_j
\end{eqnarray}
for $i,j\in\{1,2,3\}$ of the condition 
$\Lie_E(\circ)=1\cdot \circ$
(where we assume the multiplication to be as in \eqref{5.1}).
 \hfill$\Box$

\begin{remarks}\label{t5.3}
(i) The ideal $(\paa_2b_2)\subset\C\{t_2,t_3\}$
up to coordinate changes is a rich invariant.
It shows that there is a functional parameter in the
family of 3-dimensional germs of $F$-manifolds 
with $T_tM\cong Q^{(1)}$ for all $t\in M$.

\medskip
(ii) But all these germs except the one with $b_2=0$
have the unpleasant property
$\{\sqrt{\II_M},\sqrt{\II_M}\}\not\subset \sqrt{\II_M}$:
Here $\sqrt{\II_M}=(y_1-1,y_2,y_3-b_2\}$ and 
\begin{eqnarray*}
\{y_2,y_3-b_2\}=-\paa_2(b_2),
\end{eqnarray*}
and this is in $\sqrt{\II_M}$ only if $b_2=0$.

\medskip
(iii) Therefore the germ of an $F$-manifold with $b_2=0$
is the most important one of those in Theorem \ref{t5.2}.
In the case $b_2=0$, the compatibility condition \eqref{5.2}
for the coefficients of the Euler field is empty.
There $\varepsilon_1\in t_1+\C$,
$\varepsilon_2,\varepsilon_3\in\C\{t_2,t_3\}$ are
arbitrary. 

\medskip
(iv) Theorem \ref{t3.2} improved the classification 
of the Euler fields for $\NN_2$ in Theorem \ref{t3.1} (c)
by exploiting coordinate changes which do not change the
multiplication. We expect that a similar reduction 
of Euler fields to normal forms is possible for the
case $b_2=0$ in Theorem \ref{t5.2}. 
But it does not look so important. 
\end{remarks}

Next we classify the irreducible germs $(M,0)$ of $F$-manifolds
with $T_tM\cong Q^{(2)}$ for generic (or all) $t\in M$.
It is also surprisingly rich. There is also 
a functional parameter.

\begin{theorem}\label{t5.4}
The following three constructions give (up to isomorphism)
all germs $(M,0)$ of 3-dimensional $F$-manifolds with
$T_tM\cong Q^{(2)}$ for generic $t\in M$. 
The three constructions do not overlap. 
${\bf m}\subset\C\{t_2,t_3\}$ denotes the maximal ideal.

\medskip
(a Up to isomorphism, there is only one germ of a 
3-dimensional $F$-manifold with $T_tM\cong Q^{(2)}$ for
all $t\in M$. In suitable coordinates $t=(t_1,t_2,t_3)$ 
it looks as follows.
\begin{eqnarray}\label{5.6}
(M,0)=(\C^3,0),\ e=\paa_1,\ \paa_2^{\circ 2}=\paa_3,\ 
\paa_2\circ\paa_3=\paa_3^{\circ 2}=0.
\end{eqnarray}
An Euler field is a vector field of the shape
\begin{eqnarray}\label{5.7}
E&=&(t_1+c_1)\paa_1+\varepsilon_2\paa_2+(\varepsilon_{3,0}
+t_3(2\paa_2\varepsilon_2-1))\paa_3\\
&&\textup{with}\quad c_1\in\C,\ \varepsilon_2,
\varepsilon_{3,0}\in\C\{t_2\}.
\nonumber
\end{eqnarray}

\medskip 
(b) Consider an arbitrary $f\in{\bf m}-\{0\}$.
Then $(M,0)=(\C^3,0)$ with $e=\paa_1$ and with the 
multiplication given by 
\begin{eqnarray}\label{5.8}
\paa_2^{\circ 2}=f\cdot\paa_3,\ 
\paa_2\circ\paa_3=\paa_3^{\circ 2}=0
\end{eqnarray}
is an $F$-manifold with $T_tM\cong Q^{(2)}$ for generic 
$t\in M$. Here $\C\times f^{-1}(0)
=\{t\in M\,|\, T_tM\cong Q^{(1)}\}$.
The ideal $(f)\subset\{\bf m\}$ is up to coordinate changes
an invariant of the germ of an $F$-manifold.
An Euler field is a vector field of the shape
\begin{eqnarray}\label{5.9}
E&=&(t_1+c_1)\paa_1+\varepsilon_2\paa_2+\varepsilon_3\paa_3\\
&&\textup{with}\quad c_1\in\C,\ \varepsilon_2\in\C\{t_2\},
\ \varepsilon_3\in\C\{t_2,t_3\}\quad\textup{and}\nonumber\\
0&=&(\varepsilon_2\paa_2+\varepsilon_3\paa_3)(f)
+f(2\paa_2(\varepsilon_2)-\paa_3(\varepsilon_3)-1).\label{5.10}
\end{eqnarray}

\medskip
(c) Consider arbitrary $f_1,f_2\in{\bf m}$
with $\gcd(f_1,f_2)=1$ and an arbitrary 
$h\in\C\{t_2,t_3\}-\{0\}$. 
Define for $(M,0)=(\C^3,0)$ the vector field
$\sigma:=hf_2\paa_2+hf_1\paa_3$. 
Then $(M,0)$ with $e=\paa_1$ and
with the multiplication given by 
\begin{eqnarray}\label{5.11}
\paa_2^{\circ 2}=f_1^2\sigma,\ \paa_2\circ\paa_3=-f_1f_2\sigma,
\ \paa_3^{\circ 2}=f_2^2\sigma,
\end{eqnarray}
is an $F$-manifold with $T_tM\cong Q^{(2)}$ for generic 
$t\in M$. Here 
$\{t\in M\,|\, T_tM\cong Q^{(1)}\}$ is equal to 
$\C\times h^{-1}(0)$ if $h\in{\bf m}$, 
and equal to $\C\times\{0\}$ if $h(0)\neq 0$. 
Here $\paa_2\circ\sigma=\paa_3\circ\sigma=0$. 
The ideals $(f_1,f_2)\subset{\bf m}$  and 
$(h)\subset\C\{t_2,t_3\}$ are up to 
coordinate changes invariants of the germ of an $F$-manifold.
An Euler field is a vector field of the shape
\begin{eqnarray}\label{5.12}
E&=&(t_1+c_1)\paa_1+\varepsilon_2\paa_2+\varepsilon_3\paa_3\\
&&\textup{with}\quad c_1\in\C,\ \varepsilon_2,
\varepsilon_3\in\C\{t_2,t_3\}\quad\textup{and}\nonumber\\
0&=& 3h(\varepsilon_2\paa_2+\varepsilon_3\paa_3)(f_1)
+f_1(\varepsilon_2\paa_2+\varepsilon_3\paa_3)(h)\nonumber\\
&&+2f_1\paa_2(\varepsilon_2)-3f_2\paa_2(\varepsilon_3)
-f_1\paa_3(\varepsilon_3)-f_1,\label{5.13}\\
0&=& 3h(\varepsilon_2\paa_2+\varepsilon_3\paa_3)(f_2)
+f_2(\varepsilon_2\paa_2+\varepsilon_3\paa_3)(h)\nonumber\\
&&+2f_2\paa_3(\varepsilon_3)-3f_1\paa_3(\varepsilon_2)
-f_2\paa_2(\varepsilon_2)-f_2.\label{5.14}
\end{eqnarray}
\end{theorem}

{\bf Proof:}
Let $((M,0),\circ,e)$ be a germ of an $F$-manifold with
$T_tM\cong Q^{(2)}$ for generic $t\in M$. 
We choose coordinates $t=(t_1,t_2,t_3)$ with $t(0)=0$.
and use the Notations \ref{t4.1} and Lemma \ref{t4.2}.
By Lemma \ref{t4.3}, $(R_1,R_2,R_3)=0$, but
$(a_3,a_2,c_3,c_2)\neq 0$.
We are in the case \eqref{4.30} in Lemma \ref{t4.4}.

The coordinates can and will be chosen as in Remark \ref{t4.5}.
Therefore for $\paa_2\circ$ as well as for $\paa_3\circ$,
the sum of the eigenvalues is 0. As each algebra $T_tM$
is irreducible, in both cases there is only one eigenvalue.
Therefore it is 0, and $\paa_2\circ$ and $\paa_3\circ$
are nilpotent.

For generic $t$, $T_tM\cong Q^{(2)}$, and at least one of
$\paa_2|_t$ and $\paa_3|_t$ is not in the (1-dimensional)
socle of the algebra $T_tM$. Suppose $\paa_2|_t$
is not in the socle. Then $\paa_2|_t\circ\paa_2|_t\neq 0$,
but it is in the socle. Therefore the section 
$\paa_2\circ\paa_2$ is $\neq 0$, and for any $t\in M$
its value is in the socle of $T_tM$ (remark that 0 is in the
socle). Write
$\paa_2\circ\paa_2=\www f_2\paa_2+\www f_1\paa_3$
with $\www f_1,\www f_2\in\C\{t_2,t_3\}$. 

Recall that $\C\{t_2,t_3\}$ is a factorial ring
(e.g. \cite[Theorem 1.16]{GLS07}). Divide out joint 
factors of $\www f_1$ and $\www f_2$ and obtain a section
$\rho:=f_2\paa_2+f_1\paa_3$ with $\gcd(f_1,f_2)=1$.
Then for each $t\in M-\C\times\{0\}$, $\gcd(f_1,f_2)=1$
implies that $\rho|_t\neq 0$, and 
$\rho|_t$ is in the socle of $T_tM$. 
Therefore any section $\www\rho$ with $\www\rho|_t$ in the socle
for $t\in M-\C\times\{0\}$ has the shape
$\www\rho=g\cdot\rho$ with 
$g\in\OO_{M-\C\times\{0\}}=\OO_M$. 
Especially $\paa_2^{\circ 2},\paa_2\circ\paa_3,\paa_3^{\circ 2}
\in \OO_M\cdot \rho$. 

Now we consider two cases. In the 1st case, 
$(f_1,f_2)(0)\neq (0,0)$, and then we suppose first 
$f_1(0)\neq 0$, and then by multiplying $\rho$ with a unit,
we can arrange $f_1=1$. In the 2nd case $(f_1,f_2)(0)=(0,0)$. 

\medskip
{\bf 1st case,} $f_1=1$: 
We make a coordinate change $t=t(s)$ with $t_1=s_1$, 
$t_3=s_3$ and $t_2=t_2(s_2,s_3)$ such that 
$\www\paa_3 t_2(s)=f_2(t(s))$. Then
\begin{eqnarray*}
\www\paa_2&=&\www\paa_2 t_2\cdot\paa_2
+\www\paa_2 t_3\cdot\paa_3=\www\paa_2 t_2\cdot\paa_2,\\
\www\paa_3&=&\www\paa_3 t_2\cdot\paa_2
+\www\paa_3 t_3\cdot\paa_3
=f_2(t(s))\cdot\paa_2+\paa_3=\rho(t(s)).
\end{eqnarray*}
We call the new coordinates again $t$. Then $\paa_3=\rho$.
This shows \eqref{5.8} for a function $f\in\C\{t_2,t_3\}$.

In the case $f(0)\neq 0$, a coordinate change $t=t(s)$
with $t_1=s_1$, $t_2=s_2$ and $t_3=t_3(s_2,s_3)$ such that 
$\www\paa_3 t_3=f(t(s))$ exists and gives 
\begin{eqnarray*}
\www\paa_2=\paa_2+\www\paa_2 t_3\cdot\paa_3, \quad
\www\paa_3=\www\paa_3 t_3\cdot \paa_3
\quad\textup{with}\quad\www\paa_3 t_3\in\C\{s_2,s_3\}^*,\\
\www\paa_2^{\circ 2}=\paa_2^{\circ 2}
=f(t(s))\cdot\paa_3=f(t(s))(\www\paa_3 t_3)^{-1}\cdot\www\paa_3
=\www\paa_3,
\end{eqnarray*}
so we obtain \eqref{5.6}. 

In order to show that the ideal $(f)$ up to coordinate
changes is an invariant of the germ $(M,0)$ of an $F$-manifold,
we have to consider all coordinate changes which respect
the shape of \eqref{5.8}. These are coordinate changes
such that $\www\paa_3$ is a multiple by a unit of $\paa_3$,
and $\www\paa_2\circ$ is still nilpotent. Thus
$t_1=s_1$ and $t_2=t_2(s_2,s_3)$ such that $\www\paa_3t_2=0$,
so $t_2=t_2(s_2)$. 
Then $\www\paa_2t_2$ and $\www\paa_3t_3$ are units in
$\C\{t_2,t_3\}$, and 
\begin{eqnarray*}
\www\paa_2=\www\paa_2t_2\cdot\paa_2+\www\paa_2t_3\cdot\paa_3, 
\quad \www\paa_3=\www\paa_3t_3\cdot\paa_3,\\
\www\paa_2^{\circ 2}=(\www\paa_2t_2)^2(\www\paa_3t_3)^{-1}
\cdot f\cdot\www\paa_3,\quad\textup{so}\quad 
\www f=(\www\paa_2t_2)^2(\www\paa_3t_3)^{-1}\cdot f(t(s)).
\end{eqnarray*}
$f$ and $\www f$ generate the same ideal up to a coordinate
change. 

Now let $(M,0)$ be the germ of a manifold with 
the multiplication in \eqref{5.8} for some 
$f\in\C\{t_2,t_3\}-\{0\}$ on $TM$. 
We have to show that it is an $F$-manifold.
With the notations in \eqref{4.4}--\eqref{4.6}, we have
$a_3=f$ and $a_2=b_2=b_3=c_2=c_3=a_1=b_1=c_1$
and therefore $A_2=A_2^{dual}=A_3=0$ in \eqref{4.27}.
Lemma \ref{t4.3} applies and shows that $M$ is an
$F$-manifold. Also Lemma \ref{t4.4} applies,
the vanishing of $R_1,R_2,R_3$ and the nonvanishing of
$a_3=f$ show $T_tM\cong Q^{(2)}$ for generic $t\in M$.

For the shape of the Euler field
$E=\varepsilon_1\paa_1+\varepsilon_2\paa_2+\varepsilon_3\paa_3$,
one has to study the explicit version \eqref{5.5} of the
condition $\textup{Lie}_E(\circ)=1\cdot \circ$. 
The case $(i,j)=(1,1)$ gives $[e,E]=e$ and 
$\varepsilon_j\in\delta_{j1}+\C\{t_2,t_3\}$.
The cases $(i,j)\in \{(2,1),(3,1)\}$ give nothing. 
The case $(i,j)=(3,3)$ gives $\paa_3\varepsilon_1=0$.
The case $(i,j)=(2,3)$ gives this again and additionally
$\paa_2\varepsilon_1+f\paa_3\varepsilon_2=0$.
The case $(i,j)=(2,2)$ gives 
$2\paa_2\varepsilon_1-f\paa_3\varepsilon_2=0$ and \eqref{5.10}.
We obtain \eqref{5.9} and \eqref{5.10}.
The case $f=1$ specializes this to \eqref{5.7}.
The parts (a) and (b) are proved.

\medskip
{\bf 2nd case,} $f_1,f_2\in{\bf m}$: We have
\begin{eqnarray*}
\paa_2^{\circ 2}=g_1\rho,\ \paa_2\circ\paa_3=g_2\rho,\ 
\paa_3^{\circ 2}=g_3\rho\quad\textup{for some}\quad
g_1,g_2,g_3\in\C\{t_2,t_3\}.
\end{eqnarray*}
One calculates
\begin{eqnarray*}
0&=&\paa_2\circ\rho = (f_2g_1+f_1g_2)\rho,
\quad\textup{so}\quad 0=f_2g_1+f_1g_2,\\
0&=&\paa_3\circ\rho = (f_2g_2+f_1g_3)\rho,
\quad\textup{so}\quad 0=f_2g_2+f_1g_3.
\end{eqnarray*}
As $\C\{t_2,t_3\}$ is a factorial ring and
$\gcd(f_1,f_2)=1$, this implies
\begin{eqnarray*}
(g_1,g_2,g_3)=(hf_1^2,-hf_1f_2,hf_2^2)\quad\textup{for some}
\quad h\in\C\{t_2,t_3\}.
\end{eqnarray*}
We define $\sigma:=h\cdot \rho$ and obtain the multiplication 
in \eqref{5.11}. 

In order to show that the ideals $(f_1,f_2)$ and $(h)$
up to coordinate changes are invariants of the germ $M$
of an $R$-manifold, we have to consider all coordinate
changes which respect the shape of \eqref{5.11}.
The arguments above show that it is sufficient that
$\www\paa_2\circ$ and $\www\paa_3\circ$ are nilpotent.
Therefore we consider a coordinate change which satisfies
$t_1=s_1$ and $(t_2,t_3)=(t_2(s_2,s_3),t_3(s_2,s_3))$. Then
\begin{eqnarray*}
\textup{(unit)}(f_2\paa_2+f_1\paa_3)
&=&\textup{(unit)}\cdot\rho(t(s))=\www\rho
=\www f_2\www\paa_2+\www f_1\www\paa_3\\
&=& (\www f_2\www\paa_2t_2+\www f_1\www\paa_3t_2)\paa_2
+ (\www f_2\www\paa_2t_3+\www f_1\www\paa_3t_3)\paa_3.
\end{eqnarray*}
The equality $(f_2(t(s)),f_1(t(s)))=(\www f_2,\www f_1)$
of ideals follows. 
Consider the set $\{g\in\C\{t_2,t_3\}\,|\, \exists
\textup{ a vector field }X\textup{ with }[e,X]=0\textup{ and }
X\circ \textup{nilpotent and }X^{\circ 2}=g\rho\}$. 
The function $h$ is a
greatest common divisor of all functions $g$ in this set.
This property is coordinate independent.
Therefore the ideal $(h)$ up to coordinate changes is an 
invariant of the germ $(M,0)$ of an $F$-manifold. 

Now let $(M,0)$ be the germ of a manifold with the
multiplication in \eqref{5.11} for some $f_1,f_2\in{\bf m}$
with $\gcd(f_1,f_2)=1$ and for some $h\in\C\{t_2,t_3\}$.
We have to show that it is an $F$-manifold. 
One calculates immediately $\paa_2\circ\sigma=0$ and
$\paa_3\circ\sigma=0$, so the section $\sigma$
is everywhere in the socle. Because of this and 
\eqref{5.11}, $\paa_2\circ$ and $\paa_3\circ$ are nilpotent. 
Calculation and comparison with \eqref{4.19}--\eqref{4.21}
and Remark \ref{t4.5} give
\begin{eqnarray*}
\paa_2^{\circ 2}&=& hf_1^2f_2\paa_2+hf_1^3\paa_3
=\frac{1}{3}a_2\paa_2 + a_3\paa_3,\\
\paa_2\circ\paa_3&=& -hf_1f_2^2\paa_2-hf_1^2f_2\paa_3
=-\frac{1}{3}c_3\paa_2 -\frac{1}{3}a_2\paa_3,\\
\paa_2^{\circ 2}&=& hf_2^3\paa_2+hf_1f_2^2\paa_3
=c_2\paa_2 + \frac{1}{3}c_3\paa_3,
\end{eqnarray*}
so 
\begin{eqnarray}\label{5.15}
\begin{array}{c|c|c|c|c|c}
a_2 & a_3 & b_2=-\frac{1}{3}c_3& b_3=-\frac{1}{3}a_2 & c_2 & 
c_3 \\ \hline
3hf_1^2f_2 & hf_1^3 & -hf_1f_2^2 & 
-hf_1^2f_2 & hf_2^3 & 3hf_1f_2^2
\end{array}
\end{eqnarray}
Easy calculations show $A_2=A_2^{dual}=A_3=0$ for 
$A_2$, $A_2^{dual}$ and $A_3$ as in Lemma \ref{t4.3}
(or, better, in Remark \ref{t4.5}). 
Lemma \ref{t4.3} applies and shows that $M$ is an $F$-manifold.
Also Lemma \ref{t4.4} applies, the vanishing of $R_1,R_2,R_3$
and the nonvanishing of $a_2,a_3,c_2,c_3$ show 
$T_tM\cong Q^{(2)}$ for generic $t\in M$. 

For the shape of the Euler field
$E=\varepsilon_1\paa_1+\varepsilon_2\paa_2+\varepsilon_3\paa_3$,
one has to study the explicit version \eqref{5.5} of the
condition $\textup{Lie}_E(\circ)=1\cdot \circ$. 
The case $(i,j)=(1,1)$ gives $[e,E]=e$ and 
$\varepsilon_j\in_{j1}+\C\{t_2,t_3\}$.
The cases $(i,j)\in \{(2,1),(3,1)\}$ give nothing. 
The cases $(i,j)\in\{(2,2),(2,3),(3,3)\}$ give with some
tedious calculations
$\paa_2\varepsilon_1=0$, $\paa_3\varepsilon_1=0$, 
\eqref{5.13} and \eqref{5.14}. Part (c)  proved.

\hfill$\Box$

\begin{remarks}\label{t5.5}
(i) The ideal $(f)$ in part (b) and the ideals $(f_1,f_2)$ 
and $(h)$ in part (c) up to coordinate changes 
are rich invariants of the germ of an $F$-manifold. 
They show that there is a functional parameter in the
family of 3-dimensional germs of $F$-manifolds 
with $T_tM\cong Q^{(2)}$ for generic $t\in M$
if $T_0M\cong Q^{(1)}$. 
This is surprising, as part (a) says that the $F$-manifold 
is near points $t\in M$ with $T_tM\cong Q^{(2)}$ 
unique up to isomorphism. 

\medskip
(ii) If in part (c) $h$ is chosen in ${\bf m}$, then 
$h$ has a clear meaning, namely
$\C\times h^{-1}(0)=\{t\in M\,|\, T_tM\cong Q^{(1)}\}$.
The meaning of the ideal $(f_1,f_2)$ is more subtle.
It tells how the rank 1 bundle of socles of the algebras
$T_tM$ on $M-\{t\in M\,|\, T_tM\cong Q^{(1)}\}$ approaches 0.

\medskip
(iii) In all three parts of Theorem \ref{t5.4}, 
$\II_M\supset (y_2^3,y_2^2y_3,y_2y_3^2,y_3^3)$ and 
$\sqrt{\II_M}=(y_1-1,y_2,y_3)$, so here
$\{\sqrt{\II_M},\sqrt{\II_M}\}\subset \sqrt{\II_M}.$

\medskip
(iv) The $F$-manifold in part (a) is with an Euler field
$E=(t_1+c_1)\paa_1+\varepsilon_2\paa_2+\varepsilon_3\paa_3$ is 
a regular $F$-manifold if and only if $\varepsilon_2(0)\neq 0$. 
In fact, up to the choice of $c_1\in\C$,
there is only one regular germ
of a 3-dimensional and everywhere irreducible $F$-manifold
\cite[Theorem 1.3]{DH17}. A germ $(M,0)$ of an $F$-manifold
with $T_0M\cong Q^{(1)}$ has no regular Euler field,
because the socle of $Q^{(1)}$ has dimension 2. 
\end{remarks}

Next we classify the irreducible germs $(M,0)$ of $F$-manifolds
with $T_0M\cong Q^{(2)}$ and $T_tM\cong Q^{(3)}$ for 
generic $t\in M$. 

\begin{theorem}\label{t5.6}
The irreducible germs $(M,0)$ of $F$-manifolds
with $T_0M\cong Q^{(2)}$ and $T_tM\cong Q^{(3)}$ for 
generic $t\in M$ form a family with the only parameter
$p\in\Z_{\geq 2}$. For fixed $p\in\Z_{\geq 2}$, the
germ of the $F$-manifold can be given as follows:
\begin{eqnarray}
(M,0)&=&(\C^3,0)\textup{ with coordinates }t=(t_1,t_2,t_3),
e=\paa_1, \nonumber\\
\paa_2^{\circ 2}&=&\varphi^2\cdot\paa_3,\ 
\paa_2\circ\paa_3=t_2^{p-2}\varphi\cdot\paa_3,\ 
\paa_3^{\circ 2}=t_2^{2p-2}\cdot\paa_3,\label{5.16}\\
\textup{with }\varphi&:=& p+(2p-2)t_2^{p-2}t_3.\nonumber
\end{eqnarray}
The caustic is $\KK=\{t\in M\,|\, t_2=0\}=\{t\in M\,|\, 
T_tM\cong Q^{(2)}\}$. A vector field 
$E=\varepsilon_1\paa_1+\varepsilon_2\paa_2+\varepsilon_3\paa_3$
is an Euler field if and only if 
$\varepsilon_1,\varepsilon_2$ and $\varepsilon_3$ have
the following shape, here $c_1\in\C$ and
$\varepsilon_{3,0}\in\C\{t_2\}$ are arbitrary, 
\begin{eqnarray}
\varepsilon_1&=&t_1+c_1\nonumber\\
\varepsilon_2&=&t_2p^{-1}(1-t_2^{p-2}\varepsilon_{3,0}),
\label{5.17}\\
\varepsilon_3&=&\varepsilon_{3,0}+t_3p^{-1}(
2-p+(2p-2)t_2^{p-2}\varepsilon_{3,0}).\nonumber
\end{eqnarray}
\end{theorem}

The following Remarks \ref{t5.7} make the geometry of the 
$F$-manifolds in Theorem \ref{t5.6} more transparent. 
The proof of Theorem \ref{t5.6} will be given after
these remarks and will contain the proof of these remarks. 

\begin{remarks}\label{t5.7}
Let $((M,0),\circ,e)$ be one of the germs of $F$-manifolds in
Theorem \ref{t5.6}. On $M-\KK$, the bundle $TM$
of algebras decomposes into the direct sum 
$TM|_{M-\KK}=T_{A_1}\oplus T_{\NN_2}$ of bundles of 
algebras isomorphic to $\C$ respectively to 
$\C[x]/(x^2)$. Write $\sigma|_{A_1}$ respectively
$\sigma|_{\NN_2}$ for the summands of a section $\sigma$
of $TM|_{M-\KK}$ in $T_{A_1}$ respectively in $T_{\NN_2}$. Then 
\begin{eqnarray}\label{5.18}
(\paa_2|_{\NN_2})^{\circ 2}&=&0,\ \paa_3|_{\NN_2}=0\\
\paa_2|_{A_1}&=&\paa_2f\cdot e|_{A_1},
\ \paa_3|_{A_1}=\paa_3f\cdot e|_{A_1}\label{5.19}\\
\textup{with }f&=&t_2^p+t_2^{2p-2}t_3,\label{5.20}\\
\textup{so }\paa_2f&=&t_2^{p-1}\varphi,\ 
\paa_3f=t_2^{2p-2}.\nonumber
\end{eqnarray}
Also
\begin{eqnarray}\label{5.21}
\paa_2^{\circ 3}&=& \paa_2f\cdot \paa_2^{\circ 2},\\
\paa_3&=&h_2\cdot \paa_2^{\circ 2}\quad\textup{with}
\quad h_2= \varphi^{-2}.\label{5.22}
\end{eqnarray}
The Euler field has freedom in $\varepsilon_{3,0}\in\C\{t_2\}$.
This is not obvious, but it is also not surprising, as
at $t\in M-\KK$  the germ of the $F$-manifold is $A_1\NN_2$,
and the Euler fields of $\NN_2$ have a similar freedom,
see \eqref{3.4}. Though in the case of $\NN_2$, one can
normalize the Euler fields by changing the variable $t_2$,
see Theorem \ref{t3.2}. This is not possible here. 
The functional freedom in $\varepsilon_{3,0}\in\C\{t_2\}$
cannot be get rid of.
\end{remarks}

{\bf Proof of Theorem \ref{t5.6} and the Remarks \ref{t5.7}:}
Let $(M,0)$ be a germ of a 3-dimensional $F$-manifold with
$T_0M\cong Q^{(2)}$ and $T_tM\cong Q^{(3)}$ for generic
$t\in M$. Then the caustic $\KK$ is a hypersurface
and $\KK=\{t\in M\,|\, T_tM\cong Q^{(2)}\}$ and
$M-\KK=\{t\in M\,|\, T_tM\cong Q^{(3)}\}$. 
The bundle $TM|_{M-\KK}$ decomposes into 
$T_{A_1}\oplus T_{\NN_2}$ as described in the Remarks 
\ref{t5.7}, and we write $\sigma=\sigma|_{A_1}+\sigma|_{\NN_2}$
for the summands of a section $\sigma$ of $TM|_{M-\KK}$. 

Because of $T_0M\cong Q^{(2)}$, Lemma \ref{t4.6} applies.
We choose the coordinates $t$ as in this lemma, so with
\eqref{4.37}--\eqref{4.43} for suitable coefficients
$g_2,g_1,g_0,h_2,h_1,h_0\in\C\{t_2,t_3\}$. 
Let ${\bf m}\subset T_0M$ be the maximal ideal in
$T_0M$. Refining \eqref{4.37}, we can choose the coordinates
$t$ even such that
\begin{eqnarray}\label{5.23}
\C\cdot\paa_2|_0\oplus\C\cdot \paa_2^{\circ 2}|_0
={\bf m}\subset T_0M,\\
\C\cdot\paa_3|_0=\C\cdot\paa_2^{\circ 2}={\bf m}^3\subset
T_0M\nonumber
\end{eqnarray}
holds. Then $h_2(0)\neq 0$ (and also 
$g_2(0)=g_1(0)=g_0(0)=h_1(0)=h_0(0)=0$). 
In the following, five coordinate changes will be made in
order to reach the normal form in Theorem \ref{t5.6}.

The eigenvalue of $\paa_2|_{\NN_2}\circ$ is a holomorphic
function on $M-\KK$ which extends continuously 
and thus holomorphically to $M$. It can be written as
$t_1+\lambda$ with $\lambda\in\C\{t_2,t_3\}$. 
We make a special coordinate change as in \eqref{4.9},
namely we choose the new coordinates $s=(s_1,s_2,s_3)$ such that
\begin{eqnarray*}
t_2=s_2,t_3=s_3,t_1=s_1+\tau\quad\textup{with}
\quad \tau\in\C\{t_2,t_3\}=\C\{s_2,s_3\}\\
\textup{ such that }\paa_2\tau=-\lambda\in\C\{t_2,t_3\}.
\end{eqnarray*}
We obtain
\begin{eqnarray*}
\www\paa_1=\paa_1,\ \www\paa_2=\paa_2+\paa_2\tau\cdot\paa_1
=\paa_2-\lambda\cdot\paa_1,\ \www\paa_3=\paa_2+\paa_3\tau
\cdot\paa_1.
\end{eqnarray*}
We call the new coordinates again $t$ and denote also the
new coefficients again as $g_j,h_j$. Now 
$\paa_2|_{\NN_2}\circ$ is nilpotent. This implies
$g_1=g_0=0$ and $\paa_2|_{A_1}=g_2\cdot e|_{A_1}$. 

The term in square brackets in line \eqref{4.42} in Lemma
\ref{t4.6} vanishes because $(M,0)$ is an $F$-manifold.
Because of $g_1=g_0=0$ it boils down to $g_2h_{02}=0$.
But $g_2\neq 0$ because $T_tM\cong Q^{(3)}$ for generic $t$.
Therefore $\paa_2h_0=0$. We make again a special coordinate
change as in \eqref{4.9}, now with $\tau\in\C\{t_3\}$
such that $\paa_3\tau=-h_0$. Then
\begin{eqnarray*}
\www\paa_1=\paa_1,\www\paa_2=\paa_2,
\www\paa_3=\paa_3+\paa_3\tau\cdot \paa_1=\paa_3-h_0\paa_1.
\end{eqnarray*}
We call the new coordinates again $t$ and denote also the
new coefficients again as $g_j,h_j$. Now $g_1=g_0=h_0=0$. 

We make a coordinate change $t=t(s)$ with $t_1=s_1$, 
$t_3=s_3$ and $t_2=t_2(s)$ with 
$\www\paa_3t_2(s)=-h_1(t(s))$. Then
\begin{eqnarray*}
\www\paa_1=\paa_1,
\www\paa_2=\www\paa_2 t_2\cdot\paa_2,
\www\paa_3=\www\paa_3 t_2\cdot\paa_2+\www\paa_3t_3\cdot\paa_3
=\paa_3-h_1\paa_2.
\end{eqnarray*}
Here $\www\paa_2t_2$ is a unit. We call the new coordinates
again $t$ and denote also the new coefficients as $g_j,h_j$.
Now $g_1=g_0=h_0=h_1=0$ and $\paa_3=h_2\paa_2^{\circ 2}$. 
And now $\paa_3|_{\NN_2}=0$ because $\paa_2|_{\NN_2}\circ$ 
is nilpotent. 
 
The term in square brackets in line \eqref{4.41}
in Lemma \ref{t4.6} vanishes. Here it is 
$\paa_2(h_2g_2^2)-\paa_3g_2$. Therefore a function
$f\in\C\{t_2,t_3\}$ with 
\begin{eqnarray}\label{5.24}
\paa_2f=g_2\quad\textup{and}\quad\paa_3f=h_2g_2^2=h_2(\paa_2f)^2
\end{eqnarray}
exists. Here $f|_{t_3=0}$ has a zero of an order 
$p\in\Z_{\geq 2}$ because $g_2(0)=0$. 

A coordinate change $t=t(s)$ with $t_1=s_1,t_3=s_3$ and 
$t_2=t_2(s_2)\in\C\{t_2\}$ exists such that 
(after calling the new coordinates again $t$)
\begin{eqnarray*}
f|_{t_3=0}=t_2^p,\quad\textup{so}\quad
f=t_2^p+\sum_{k\geq 1}f_kt_3^k\textup{ with }
f_k\in\C\{t_2\}.
\end{eqnarray*}
The equation $\paa_3f=h_2(\paa_2f)^2$ and $h_2(0)\neq 0$ and
$p\geq 2$ imply inductively $f_k\in t_2^{2p-2}\C\{t_2\}$
for all $k\geq 1$ and especially 
$f_1\in t_2^{2p-2}\cdot\C\{t_2\}^*$. Therefore
\begin{eqnarray*}
f=t_2^p\Bigl(1+t_2^{p-2}t_3\cdot 
(\textup{a unit in }\C\{t_2\}\Bigr).
\end{eqnarray*}
We can and will change the coordinate $t_3$ such that 
$f=t_2^p(1+t_2^{p-2}t_3)$. 
Then $g_2=\paa_2f=t_2^{p-1}(p+(2p-2)t_2^{p-2}t_3)
=t_2^{p-1}\varphi$ and $\paa_3f=t_2^{2p-2}$. 
Now $h_2$ is determined by $\paa_3f=h_2(\paa_2f)^2$
and is $h_2=\varphi^{-2}$. 

Now all statements in the Remarks \ref{t5.7} 
except those on the Euler field are shown. 
The terms in the square brackets in \eqref{4.41}--\eqref{4.43}
vanish. Therefore we really have an $F$-manifold. 
The multiplication is as in \eqref{5.16}, because
\begin{eqnarray*}
\paa_2^{\circ 2}&=&h_2^{-1}\paa_3=\varphi^2\paa_3,\\
\paa_2^{\circ 3}&=&g_2\paa_2^{\circ 2}
=t_2^{p-1}\varphi\paa_2^{\circ 2},\\
\paa_2\circ\paa_3&=& h_2\paa_2^{\circ 3}=h_2g_2\paa_2^{\circ 2}
=g_2\paa_3=t_2^{p-1}\varphi\paa_3,\\
\paa_3^{\circ 2}&=& h_2^2\paa_2^{\circ 4}=
h_2^2g_2^2\paa_2^{\circ 2}=h_2g_2^2\paa_3=t_2^{2p-2}\paa_3.
\end{eqnarray*}

It remains to show the shape \eqref{5.17} of the Euler field
$E=\varepsilon_1\paa_1+\varepsilon_2\paa_2+\varepsilon_3\paa_3$.
One has to study the explicit version \eqref{5.5}
of the condition $\textup{Lie}_E(\circ)=1\cdot \circ$.
The case $(i,j)=(1,1)$ just gives $[e,E]=e$ and thus 
$\varepsilon_j\in\delta_{1j}+\C\{t_2,t_3\}$. 
The cases $(i,j)\in\{(2,1),(3,1)\}$ give nothing.
The case $(i,j)=(3,3)$ leads to 
$\paa_3(\varepsilon_1)=0$, $\paa_3\varepsilon_2=0$
and $\varepsilon_2=(2p-2)^{-1}t_2(1-\paa_3\varepsilon_3)$.
The cases $(i,j)\in\{(2,2),(2,3)\}$ lead to 
$\paa_2\varepsilon_1=0$ and to equations which allow to
relate $\varepsilon_{3,0}$ and $\varepsilon_{3,1}\in\C\{t_2\}$
in $\varepsilon_3=\varepsilon_{3,0}+t_3\varepsilon_{3,1}$. 
At the end one obtains \eqref{5.17}.
We leave the details to the reader.
\hfill$\Box$ 

\bigskip
We do not have a classification of the irreducible
germs $(M,0)$ of $F$-manifolds with
$T_0M\cong Q^{(1)}$ and $T_tM\cong Q^{(3)}$ for generic
$t\in M$. The family of examples in the next lemma shows 
that such germs exist.

\begin{lemma}\label{t5.8}
Fix a number $p\in\Z_{\geq 2}$. 
The manifold $M=\C^3$ with coordinates $t=(t_1,t_2,t_3)$
and with the multiplication on $TM$ given by $e=\paa_1$ and
\begin{eqnarray}\label{5.25}
\paa_2^{\circ 2}=pt_2^{p-1}\cdot\paa_2,\ 
\paa_2\circ\paa_3=0,\ \paa_3^{\circ 2}=0
\end{eqnarray}
is an $F$-manifold with caustic
\begin{eqnarray}
\KK&=&\{t\in M\,|\, t_2=0\}=\{t\in M\,|\, T_tM\cong Q^{(1)}\}
\textup{ and}\nonumber\\
M-\KK&=&\{t\in M\,|\, T_tM\cong Q^{(3)}\}.\label{5.26}
\end{eqnarray}
A vector field $E$ is an Euler field if and only if
\begin{eqnarray}\label{5.27}
E=(t_1+c_1)\paa_1+\frac{1}{p}t_2\paa_2+\varepsilon_3\paa_3
\quad\textup{with }c_1\in\C,\ \varepsilon_3\in\C\{t_3\}.
\end{eqnarray}
\end{lemma}

{\bf Proof:}
In the Notations \ref{t4.1}, $b_2=b_3=a_3=c_2=c_3=0$,
$a_2=t_2^{p-1}$. Therefore $A_2=A_2^{dual}=A_3=0$ 
in Lemma \ref{t4.4}, and we have an $F$-manifold.
The statement on $\KK$ is clear. The analytic spectrum is
\begin{eqnarray}\label{5.28}
L_M=\{(y,t)\in T^*M\,|\, y_1=1,y_2(y_2-pt_2^{p-1})=y_2y_3=y_3^2
=0\}.
\end{eqnarray}
The set which underlies $L_M$ is 
$L_M^{red}=\{(y,t)\in T^*M\,|\, (y_1,y_2,y_3)=(1,0,0)\}
\cup\{(y,t)\in T^*M\,|\, (y_1,y_2,y_3)=(1,pt_2^{p-1},0)\}$,
so it has two components which meet over $\KK$. 
Therefore $T_tM\cong Q^{(3)}$ for $t\in M-\KK$. 
For the proof of \eqref{5.27}, one has to study the
explicit version \eqref{5.5} of the condition
$\textup{Lie}_E(\circ)=1\cdot \circ$. 
We leave the details to the reader. 
\hfill$\Box$

\section{Examples of 3-dimensional generically semisimple 
$F$-manifolds}\label{c6}
\setcounter{equation}{0}

A partial classification of irreducible germs $(M,0)$ of 
3-dimensional generically semisimple $F$-manifolds 
was undertaken in \cite[ch. 5.5]{He02}, 
there in the Theorems 5.29 and 5.30.
Theorem 5.29 in \cite{He02} gave basic facts on both cases,
the case $T_{0}M\cong Q^{(2)}$ and the case 
$T_{0}M\cong Q^{(1)}$.
Theorem 5.30 classified completely those germs 
where $T_{0}M\cong Q^{(2)}$ and where the germ
$(L_M,\lambda)$ of the analytic spectrum has 3 components.

Below we first describe in the Remarks \ref{t6.1} the
strategy of the classification results in this section
and the next section.  
The Examples \ref{t6.2} rewrite the three distinguished 
$F$-manifolds $A_3$, $B_3$ and $H_3$. 
Theorem \ref{t6.3} is Theorem 5.29 from \cite{He02}. 
Lemma \ref{t6.4} and Lemma \ref{t6.5} give examples $(M,0)$
of generically semisimple $F$-manifolds with 
$T_0M\cong Q^{(1)}$. We do not have a classification of
all such germs.  Lemma \ref{t6.4} is Theorem 5.32 from
\cite{He02}.

\begin{remarks}\label{t6.1}
(i) Let $(M,0)$ be an irreducible germ of a 3-dimensional
generically semisimple $F$-manifold with analytic spectrum 
$(L_M,\lambda)$.
Here and in the following, we choose 
coordinates $t=(t_1,t_2,t_3)$ on $(M,0)$
such that $(M,0)\cong(\C^3,0)$ and $e=\paa_1$.
Then $(y_1,y_2,y_3)$ are the fiber coordinates on $T^*M$ 
which correspond to $\paa_1,\paa_2,\paa_3$.
And $\alpha=\sum_{i=1}^3y_i\ddd t_i$ is the canonical
1-form. In the following, $M$ denotes a suitable
(small) representative of $(M,0)$. 

It turns out that often the best way to arrive at
a normal form for $((M,0),\circ,e)$ is to control 
the function $F:(L_M,\lambda)\to(\C,0)$ from
Theorem \ref{t2.15} (b). It is holomorphic on
$L_M^{reg}$ and continuous on $L$, and it satisfies
\begin{eqnarray}\label{6.1}
\ddd F|_{L_M^{reg}}=\alpha|_{L_M^{reg}}.
\end{eqnarray}
We consider it as a 3-valued holomorphic function
on $M$ which is branched precisely over the caustic
$\KK\subset M$. Locally on $M-\KK$ it splits
into three holomorphic functions
$F^{(1)},F^{(2)},F^{(3)}$. We will use this notation
without specifying a simply connected region in $M-\KK$. 
This is imprecise, but not in a harmful way. 
With this notation, $F$
determines $L_M$ as follows (this rewrites \eqref{6.1}),
locally on $M-\KK$,
\begin{eqnarray}\label{6.2}
L_M=\bigcup_{j=1}^3\{(y,t)\in T^*M\,|\, 
y_i=\paa_i F^{(j)}
\textup{ for }i\in\{1,2,3\}\}.
\end{eqnarray}
Let $M^{(r)}$ be a suitable neighborhood of $0$ in the
$(t_2,t_3)$-plane $\C^2$. It can be identified with the
set of $e$-orbits of $M$. The condition $\Lie_e(\circ)=0$
implies that the multiplication, the caustic
$\KK$ and the analytic spectrum $L_M$ are invariant 
under the flow of $e$.
As $e=\paa_1$, $\KK$ induces a hypersurface 
$\KK^{(r)}\subset M^{(r)}$, and $F=t_1+f$
where $f$ is a 3-valued holomorphic function on $M^{(r)}$
which is branched along $\KK^{(r)}$. 
Locally on $M^{(r)}-\KK^{(r)}$, $f$ splits into
three holomorphic functions $f^{(1)},f^{(2)},f^{(3)}$.
The coefficients of the polynomials
$\prod_{j=1}^3(x-f^{(j)})$ and 
$\prod_{j=1}^3(x-\paa_2 f^{(j)})$ and
$\prod_{j=1}^3(x-\paa_3 f^{(j)})$ are 
univalued holomorphic functions on $M^{(r)}$, i.e.
they are in $\C\{t_2,t_3\}$. The last two polynomials
are the characteristic polynomials of
$\paa_2\circ$ and $\paa_3\circ$.

\medskip
(ii) The Euler field 
$E=\varepsilon_1\paa_1+\varepsilon_2\paa_2
+\varepsilon_3\paa_3$ on $M-\KK$, which corresponds to $F$ by 
Theorem \ref{t2.15}, is given by $F^{(j)}=E(F^{(j)})$, i.e. by
\begin{eqnarray}\label{6.3}
\varepsilon_1=t_1,\ 
f^{(j)}= \varepsilon_2\cdot\paa_2 f^{(j)}
+\varepsilon_3\cdot \paa_3 f^{(j)}.
\end{eqnarray}
$\varepsilon_2$ and $\varepsilon_3$ depend only on $(t_2,t_3)$,
but often are meromorphic along $\KK^{(r)}$. 
If they are in $\C\{t_2,t_3\}$, then $E$ extends from
$M-\KK$ to $M$. 

\medskip
(iii) Now consider the case $T_0M\cong Q^{(2)}$.
Denote by ${\bf m}\subset T_0M$ the maximal ideal in $T_0M$.
We can and will choose the coordinates $t$ such that
\begin{eqnarray}\label{6.4}
\C\cdot \paa_2|_0\oplus \C\cdot \paa_2^{\circ 2}|_0
={\bf m}\subset T_0M,\\
\C\cdot \paa_3|_0=\C\cdot\paa_2^{\circ 2}|_0 ={\bf m}^2
\subset T_0M.\nonumber
\end{eqnarray}
Then
\begin{eqnarray}\label{6.5}
\paa_3&=& h_2\cdot\paa_2^{\circ 2}+h_1\cdot\paa_2+h_0\cdot\paa_1
\end{eqnarray}
for suitable coefficients $h_2,h_1,h_0\in\C\{t_2,t_3\}$ with 
$h_2(0)\neq 0$, $h_1(0)=h_0(0)=0$. These coefficients 
are determined by
\begin{eqnarray}\label{6.6}
\paa_3 f^{(j)}
= h_2\cdot (\paa_2 f^{(j)})^2 
+h_1\cdot \paa_2 f^{(j)}+h_0.
\end{eqnarray}
Also write
$\prod_{j=1}^3(x-\paa_2 f^{(j)})
=x^3-g_2x^2-g_1x-g_0$ with $g_2,g_1,g_0\in\C\{t_2,t_3\}$. 
Then
\begin{eqnarray}
L_M=\{(y,t)\in T^*M&|& y_1=1, 
y_2^3=g_2y_2^2+g_1y_2+g_0,\nonumber\\ 
&&y_3=h_2y_2^2+h_1y_2+h_0\}.\label{6.7}
\end{eqnarray}
\end{remarks}

\begin{examples}\label{t6.2}
Here the $F$-manifolds $A_3,B_3,H_3$ from Theorem 5.22 (i) 
in \cite{He02} will be rewritten with the notions from 
the Remarks \ref{t6.1}. They arise as complex orbit spaces
of the corresponding Coxeter groups. 
Their discriminants had been studied especially by 
O.P. Shcherbak \cite{Sh88}, and their Lagrange maps
(which correspond to the $F$-manifold structures
by Theorem 3.16 in \cite{He02}) had been studied by
Givental \cite{Gi88}. 

They are simple $F$-manifolds with Euler fields with
positive weights. Their germs $(M,0)$ at 0 are the only simple
3-dimensional germs of $F$-manifolds with $T_0M\cong Q^{(2)}$,
see Theorem \ref{t6.3} (b). 

We use the notations from the Remarks \ref{t6.1}.
Though here we have $F$-manifolds $M=\C^3$, not just germs. 
The following table gives for each of the three cases
the following data: 
\begin{list}{}{}
\item[(i)] 
a 3-valued function $\xi$ on $\C^2=M^{(r)}$ (with coordinates 
$(t_2,t_3)$) which is branched along $\KK^{(r)}$.
It is given by the equation of degree 3 which it satisfies.
The equation is denoted $\xi^3-g_2\xi^2-g_2\xi-g_0=0$
with $g_2,g_1,g_0\in\C[t_2,t_3]$. 
\item[(ii)]
A weight system $(w_1,w_2,w_3)\in\Q_{>0}^3$.
\item[(iii)]
The components of $\KK^{(r)}$, and which germs of an 
$F$-manifold are at generic points of each component.
\end{list}
\begin{eqnarray*}
\begin{array}{c|c|c|c}
 & A_3 & B_3 & H_3 \\ \hline
\xi & \xi^3+2\xi t_3+t_2 & \xi(\xi^2+2\xi t_3+t_2) &
\xi^3-(2\xi t_3+t_2)^2\\ \hline  
(w_1,w_2,w_3) & (1,\frac{3}{4},\frac{1}{2}) & 
(1,\frac{2}{3},\frac{1}{3}) & (1,\frac{3}{5},\frac{1}{5})\\ 
\hline 
\KK^{(r)}:\ A_2A_1 & 27t_2^2+32t_3^3=0: & 
t_2-t_3^2: & 27t_2+32t_3^3: \\
\KK^{(r)}\textup{ 2nd comp.} & - & t_2=0:\ I_2(4)A_1 
& t_2=0:\ I_2(5)A_1
\end{array}
\end{eqnarray*}
The 3-valued function $f$ on $M^{(r)}$ with $F=t_1+f$ is
\begin{eqnarray}\label{6.8}
f&=& w_2\xi t_2 + w_3\xi^2 t_3.
\end{eqnarray}
The following identities are crucial. They will be proved
below.
\begin{eqnarray}\label{6.9}
\paa_2 f =\xi,\quad \paa_3 f=\xi^2.
\end{eqnarray}
Because of them, the Euler field is 
$E=t_1\paa_2+w_2t_2\paa_2+w_3t_3\paa_3$, and the 
multiplication and the analytic spectrum
are given as follows,
\begin{eqnarray}\label{6.10}
L_M=\{(y,t)\in T^*M\,|\, y_1=1,y_2^3=g_2y_2^2+g_1y+g_0,
y_3=y_2^2\}.
\end{eqnarray}
One sees $L\cong\C^2\times C$, where $C$ is a plane curve,
and $C$ is smooth in the case $A_3$, $C$ has two smooth 
components which intersect transversely in the case $B_3$,
and $C$ has one ordinary cusp in the case $H_3$. 
We will prove now \eqref{6.9} and the claims on the caustic.

\medskip
{\bf Proof of \eqref{6.9}:}
It is equivalent to $\ddd f = \xi\ddd t_2+\xi^2\ddd t_3$.
And this is equivalent to the claim that $L_M^{reg}$ is 
Lagrange, i.e. $\alpha|_{L_M^{reg}}$ is closed. 
And then $\int_{L_M}\alpha=t_1+f$ as a 3-valued function on $M$.
In all three cases
\begin{eqnarray*}
\ddd f = w_2\xi\ddd t_2+w_3\xi^2\ddd t_3 + 
(w_2t_2+2w_3\xi t_3)\ddd\xi,
\end{eqnarray*}
so in all three cases
\begin{eqnarray}\label{6.11}
(w_2t_2+2w_3\xi t_3)\ddd\xi 
= (1-w_2)\xi\ddd t_2+(1-w_3)\xi^2\ddd t_3
\end{eqnarray} 
has to be shown. 

The case $A_3$: Use $0=\xi^3+2\xi t_3+t_2$ to calculate
\begin{eqnarray*}
0&=&\xi\ddd (\xi^3+2\xi t_3+t_2)
= (3\xi^3+2\xi t_3)\ddd \xi + 2\xi^2\ddd t_3+\xi\ddd t_2\\
&=& (-4\xi t_3-3t_2)\ddd \xi + 2\xi^2\ddd t_3+\xi\ddd t_2,
\end{eqnarray*}
which shows \eqref{6.11}.

The case $B_3$: Use $0=\xi(\xi^2+2\xi t_3+t_2)$ 
to calculate
\begin{eqnarray*}
0&=&\xi\ddd (\xi(\xi^2+2\xi t_3+t_2))
= (3\xi^3+4\xi^2 t_3+\xi t_2)\ddd \xi 
+ 2\xi^3\ddd t_3+\xi^2\ddd t_2\\
&=& (-2\xi^2 t_3-2\xi t_2)\ddd \xi
+2\xi^3\ddd t_3+\xi^2\ddd t_2\\
&=&2\xi\cdot \eta\quad\textup{with}\quad
\eta:=(-\xi t_3-t_2)\ddd\xi + \xi^2\ddd t_3+
\frac{1}{2}\xi\ddd t_2.
\end{eqnarray*}
For \eqref{6.11}, we need $\eta=0$.
Here $\xi$ consists of one holomorphic function $\xi^{(1)}=0$
and a 2-valued function $\xi^{(2\&3)}\neq 0$. 
For $\xi^{(1)}$, $\eta=0$ is trivial.
For $\xi^{(2\&3)}$, $\eta=0$ follows from 
$0=2\xi^{(2\&3)}\cdot \eta$, 
as then we may divide by $2\xi^{(2\&3)}$.

The case $H_3$: Use $0=\xi^3-(2\xi t_3+t_2)^2$ 
to calculate
\begin{eqnarray*}
0&=& \xi\ddd (\xi^3-(2\xi t_3+t_2)^2)
= 3\xi^3\ddd \xi -2\xi(2\xi t_3+t_2)\ddd (2\xi t_3+t_2)\\
&=& (2\xi t_3+t_2)\cdot\Bigl(3(2\xi t_3+t_2)\ddd \xi 
-2\xi \ddd (2\xi t_3+t_2)\Bigr),\quad\textup{thus}\\
0&=& 3(2\xi t_3+t_2)\ddd \xi - 2\xi \ddd (2\xi t_3+t_2)\\
&=& (2\xi t_3+3t_2)\ddd \xi - 4\xi^2\ddd t_3-2\xi\ddd t_2,
\end{eqnarray*}
which shows \eqref{6.11}. \hfill$\Box$

\medskip
{\bf Proof of the statements on the caustic:}

The case $A_3$: The discriminant of $x^3+2t_3x+t_2$ is
$4(2t_3)^3+27t_2^2=32t_3^3+27t_2^2$. Over generic points 
of the caustic (all except $(t_1,0,0)$ for $t_1\in\C$)
the multigerm of $L_M$ has 2 smooth components, so there
we have the germ $A_2A_1$. 

The case $B_3$: The two components of $L_M$ meet over points
with $t_2=0$. Over generic points of this component of $\KK$,
we have the germ $I_2(4)A_1$. The discriminant of
$x^2+2t_3x+t_2$ is $(2t_3)^2-4t_2=4(t_3^2-t_2)$. Over generic
components of this component, the multigerm of $L_M$ has
2 smooth components, so there we have the germ $A_2A_1$. 

The case $H_3$: The discriminant of $x^3-(2t_3x+t_2)^2$
is $t_2^3(32t_3^3+27t_2)$. The cusp surface of $L_M$
lies over the component with $t_2=0$ of $\KK$. 
So there we have the germ $I_2(5)A_1$. Over generic points of
the component with $32t_3^3+27t_2=0$, the multigerm $L_M$
has 2 smooth components, so there we have the germ $A_2A_1$.
\hfill$\Box$ 
\end{examples}

The following theorem is Theorem 5.29 from \cite{He02}.
It gives basic facts on irreducible germs of 3-dimensional 
$F$-manifolds with generically semisimple multiplication.

\begin{theorem}\label{t6.3}\cite[Theorem 5.29]{He02}
Let $(M,0)$ be an irreducible germ of a 3-dimensional
generically semisimple $F$-manifold
with analytic spectrum $(L_M,\lambda)\subset T^{*}M$.

\medskip
(a) Suppose $T_{0}M\cong Q^{(2)}$. 
Then $(L_M,\lambda)$ has embedding dimension 3 or 4 and 
$(L_M,\lambda)\cong (\C^2,0)\times (C,0)$ for a
plane curve $(C,0)\subset (\C^2,0)$ with $\mult (C,0)\leq 3$.
The Euler field $E_0$ from Theorem \ref{t2.15} (c) on 
$M-\KK$ extends holomorphically to $M$ if and only if 
$(C,0)$ is quasihomogeneous.

\medskip
(b) Suppose $T_{0}M\cong Q^{(2)}$ and 
$(L_M,\lambda)\cong (\C^2,0)\times (C,0)$ 
with $\mult (C,0)<3$. 
Then $((M,0),\circ,e)$ is one of the germs $A_3, B_3, H_3$.

\medskip
(c) Suppose $T_{0}M\cong Q^{(2)}$ and 
$(L_M,\lambda)\cong (\C^2,0)\times (C,0)$ with $\mult(C,0)=3$. 
Then the caustic $\KK$ is a smooth surface and coincides 
with the $\mu$-constant stratum. That means, 
$T_qM\cong Q^{(2)}$ for each $q\in \KK$.
The modality is $\mmod_\mu (M,0)=1$ (the maximal possible)
(recall Definition \ref{t2.9} (a)).

\medskip
(d) Suppose $T_{0}M\cong Q^{(1)}$. Then $(L_M,\lambda)$ 
has embedding dimension 5 and 
$(L_M,\lambda)\cong (\C,0)\times (L^{(r)},0)$.
Here $(L^{(r)},0)$ is a Lagrange surface with embedding 
dimension 4. Its ring $\OO_{L^{(r)},0}$ is a 
Cohen-Macaulay ring, but not a Gorenstein ring.
\end{theorem}

{\bf Sketch of the proof:}
(a) One chooses the coordinates $(t_1,t_2,t_3)$ as in
\eqref{6.4}. Then $L_M$ is as in \eqref{6.7}.
Because of the equations $y_1=1$ and 
$y_3=\sum_{i=0}^2h_iy_2^i$, $(L_M,\lambda)$
has embedding dimension $\leq 4$. Theorem \ref{t2.14}
applies and gives $(L_M,\lambda)\cong(\C^2,0)\times (C,0)$
for a plane curve $C$. The germ $(C,0)$ has multiplicity 
$\leq 3$ because the projection $\pi_L:L_M\to M$ is a branched
covering of degree 3. 
The Euler field $E_0$ from Theorem \ref{t2.15} (c)
on $M-\KK$ extends to $M$ if and only if $(C,0)$
is quasihomogeneous because of Theorem \ref{t2.16} (e) 
and Theorem \ref{t2.15} (b)+(c).

(b) $\mult(C,0)\leq 2$ means that $(C,0)$ is either smooth
or a double point or a cusp. In the first two cases,
one can apply the correspondence between $F$-manifolds
and hypersurface or boundary singularities 
\cite[Theorem 5.6 and Theorem 5.14]{He02}
and the fact that $A_3$, $B_3$ and $C_3$ are the only
hypersurface or boundary singularities with Milnor
number 3. In the case of a cusp, results of Givental
\cite{Gi88} are used, see the proof of 
\cite[Theorem 5.29]{He02}.

(c) $(C,0)$ has multiplicity 3. And the projection 
$\pi_L:L_M\to M$ is a branched covering of degree 3. 
Together these facts imply that $\pi_L:L_M\to M$ is precisely
branched at the points of $L_M$ which correspond to
$(\C^2,0)\times\{0\}$ in 
$(\C^2,0)\times (C,0)\cong (L_M,\lambda)$. 
This implies all the statements.

(d) If the embedding dimension of $(L_M,\lambda)$ were $\leq 4$,
then by Theorem \ref{t2.14} 
$(L_M,\lambda)\cong(\C^2)\times (C,0)$ for a plane curve,
so then $(L_M,\lambda)$ were a complete intersection,
and thus $T_0M\cong Q^{(2)}$, a contradiction. 
Therefore the embedding dimension is 5 and by Theorem
\ref{t2.14} $(L_M,\lambda)\cong (\C,0)\times (L^{(r)},0)$
where $(L^{(r)},0)$ has embedding dimension 4. 
The ring $\OO_{L^{(r)},0}$ is Cohen-Macaulay because
$\pi^{(r)}:L^{(r)}\to M^{(r)}$ is finite and flat. It is not 
Gorenstein, because $T_0M\cong Q^{(1)}$ is not Gorenstein.
\hfill$\Box$

\bigskip
The classification of germs $(M,0)$ of 3-dim
generically semisimple $F$-manifolds with
$T_0M\cong Q^{(1)}$ is not treated in this paper.
Only one family of examples from \cite{He02}
and one other interesting example will be given now.

\begin{lemma}\label{t6.4}\cite[Theorem 5.32]{He02}
Fix two numbers $p_2,p_3\in\Z_{\geq 2}$. 
The manifold $M=\C^3$ with coordinates $t=(t_1,t_2,t_3)$
and with the multiplication on $TM$ given by $e=\paa_1$ and
\begin{eqnarray}\label{6.12}
\paa_2^{\circ 2}=p_2t_2^{p_2-1}\cdot\paa_2,\ 
\paa_2\circ\paa_3=0,\ \paa_3^{\circ 2}=p_3t_3^{p_3-1}\paa_3
\end{eqnarray}
is a simple (and thus generically semisimple) $F$-manifold
with $T_0M\cong Q^{(1)}$. 
Its caustic $\KK$ has two components 
$\KK^{(1)}=\{t\in M\,|\, t_2=0\}$ and 
$\KK^{(2)}=\{t\in M\,|\, t_3=0\}$. 
The germ $(M,t)$ is of type $A_1I_2(2p_2)$
for $t\in\KK^{(1)}-\C\times\{0\}$ and of type
$A_1I_2(2p_3)$ for $t\in \KK^{(2)}-\C\times\{0\}$. 
A vector field $E$ is an Euler field if and only if
\begin{eqnarray}\label{6.13}
E=(t_1+c_1)\paa_1+\frac{1}{p_2}t_2\paa_2+\frac{1}{p_3}t_3\paa_3
\quad\textup{with }c_1\in\C.
\end{eqnarray}
\end{lemma}

{\bf Proof:} The analytic spectrum is
\begin{eqnarray}\label{6.14}
L_M=\{(y,t)\in T^*M&|& y_1=1,y_2(y_2-p_2t_2^{p_2-1})=y_2y_3=0,\\
&&y_3(y_3-p_3t_3^{p_3-1})=0\}.\nonumber
\end{eqnarray}
The set which underlies $L_M$ has three components
$L^{(1)},L^{(2)},L^{(3)}$ with
\begin{eqnarray*}
L^{(1)}&=&\{(y,t)\in T^*M\,|\, (y_1,y_2,y_3)=(1,0,0)\}\\
L^{(2)}&=&\{(y,t)\in T^*M\,|\, 
(y_1,y_2,y_3)=(1,p_2t_2^{p_2-1},0)\}\\
L^{(3)}&=&\{(y,t)\in T^*M\,|\,
(y_1,y_2,y_3)=(1,0,p_3t_3^{p_3-1})\}.
\end{eqnarray*}
The functions 
$f^{(1)}:=0, f^{(2)}:=t_2^{p_2}, f^{(3)}:=t_3^{p_3}$ on $M$ 
satisfy 
\begin{eqnarray*}
f^{(j)}=\frac{1}{p_2}t_2\paa_2f^{(j)}+\frac{1}{p_3}t_3
\paa_3f^{(j)}.
\end{eqnarray*}
If one lifts $f^{(j)}$ to $L^{(j)}$, the resulting function
on $L_M$ is $\frac{1}{p_2}t_2y_2+\frac{1}{p_3}t_3y_3$, so
a holomorphic function on $L_M$.
$F:=t_1+f$ satisfies all properties in the Remarks \ref{t6.1}
(i)+(ii). Therefore $(M,\circ,e)$ is an $F$-manifold,
and the Euler field is as claimed.

$L^{(2)}$ and $L^{(3)}$ meet only over $t=0$,
$L^{(1)}$ and $L^{(2)}$ meet over $\KK^{(1)}$, 
and $L^{(1)}$ and $L^{(3)}$ meet over $\KK^{(2)}$.
From their intersection multiplicities or from the coefficients
of the Euler field, one concludes that the germs of 
$F$-manifolds have the types $A_1I_2(2p_2)$ respectively
$A_2I_2(2p_3)$ at points of $\KK^{(1)}$ respectively
$\KK^{(2)}$ not equal to $\C\times\{0\}\subset M$. 
The stratification of $M$ by the types of the germs
of the $F$-manifolds shows that the $F$-manifold is simple. 
\hfill$\Box$ 

\bigskip

Als the following example is a simple $F$-manifold
$M=\C^3$ with $T_0M\cong Q^{(1)}$. But its analytic spectrum
is irreducible,  and it is singular only in codimension 2.

\begin{lemma}\label{t6.5}
The manifold $M=\C^3$ with coordinates $t=(t_1,t_2,t_3)$
and with the multiplication on $TM$ given by $e=\paa_1$ and
\begin{eqnarray}
(\paa_2-\frac{1}{2}t_3\paa_1)^{\circ 2}=\frac{9}{4}t_2^2\paa_1
-\frac{3}{2}t_3(\paa_2-\frac{1}{2}t_3\paa_1)-
\frac{3}{2}t_2(\paa_3+\frac{1}{2}t_2\paa_1),\nonumber\\
(\paa_2-\frac{1}{2}t_3\paa_1)\circ 
(\paa_3+\frac{1}{2}t_2\paa_1)=\frac{3}{4}t_2t_3,\label{6.15}\\ 
(\paa_3+\frac{1}{2}t_2\paa_1)^{\circ 2}=-\frac{3}{4}t_3^2\paa_1
-\frac{1}{2}t_3(\paa_2-\frac{1}{2}t_3\paa_1)
+\frac{3}{2}t_2(\paa_3+\frac{1}{2}t_2\paa_1),\nonumber
\end{eqnarray}
is a simple (and thus generically semisimple) $F$-manifold
with $T_0M\cong Q^{(1)}$. 
The Lagrange surface $L^{(r)}$ in $T^*M^{(r)}$ with
$L_M\cong\C\times L^{(r)}$ is smooth outside $0$.
The caustic $\KK$ has 4 components.
The corresponding 4 components of $\KK^{(r)}\subset M^{(r)}$
are the 4 lines through 0 which are together given by
\begin{eqnarray}\label{6.16}
0=t_3^4+6t_2^2t_3^2-3t_2^4.
\end{eqnarray}
The germ $(M,t)$ is of type $A_1A_2$ 
for $t\in \KK-\C\times\{0\}$. 
A vector field $E$ is an Euler field if and only if
\begin{eqnarray}\label{6.17}
E=(t_1+c_1)\paa_1+\frac{1}{2}t_2\paa_2+\frac{1}{2}t_3\paa_3
\quad\textup{with }c_1\in\C.
\end{eqnarray}
\end{lemma}

{\bf Proof:} In the Notations \ref{t4.1} in 
\eqref{4.4}--\eqref{4.6}, 
\begin{eqnarray}\label{6.18}
(a_2,a_3,b_2,b_3,c_2,c_3)=(\frac{-3}{2}t_3,\frac{-3}{2}t_2,
\frac{-1}{2}t_2,\frac{1}{2}t_3,\frac{-1}{2}t_3,\frac{3}{2}t_2),
\\
(a_1,b_1,c_1)=(\frac{9}{4}t_2^2,\frac{3}{4}t_2t_3,
\frac{-3}{4}t_3^2)=(-a_3c_3,a_3c_2,-a_2c_2).\label{6.19}
\end{eqnarray}
Therefore the multiplication is associative. 
One checks easily that $A_2,A_2^{dual}$ and $A_3$ in 
Lemma \ref{t4.4} vanish. Therefore $(M,\circ,e)$ is an
$F$-manifold. The function $9R_3^2-4R_1R_2$ in Lemma
\ref{t4.3} is here 
$9R_3^2-4R_1R_2=\frac{9}{4}(t_3^4+6t_2^2t_3^2-3t_2^4)$.
Therefore the $F$-manifold is generically semisimple,
and the caustic is given by \eqref{6.16}. 
One checks also easily that the explicit version
\eqref{5.5} for $\textup{Lie}_E(\circ)=1\cdot\circ$
is satisfied for $E$ as in \eqref{6.17}.
Therefore $E$ is one Euler field. 
By Theorem \ref{t2.15} (b)+(c) and the irreducibility of
the germ $(M,0)$ the vector fields $E+c_1e$ for $c_1\in\C$
are the only Euler fields on $M$.

The smoothness of $L^{(r)}$ outside $0$ would
imply together with the classification of the 2-dimensional 
germs of $F$-manifolds that the germ $(M,t)$ for 
$t\in \KK-\C\times\{0\}$ is of type $A_1A_2$. 
It remains to show that $L^{(r)}$ is smooth. 
For this we reveal how the $F$-manifold
was constructed. Consider the coordinate change on 
$T^*M^{(r)}$ with new coordinates $x=(x_1,x_2,x_3,x_4)$,
\begin{eqnarray*}
t_2=x_1+x_3,\ t_3=x_2-x_4,\ y_2=x_2+x_4,\ y_3=-x_1+x_3.
\end{eqnarray*}
In the new coordinates, the functions $Y_{22},Y_{23}$
and $Y_{33}$ in Lemma \ref{t4.4}, which define $L^{(r)}$,
become
\begin{eqnarray*}
Y_{22}&=& (x_2^2-3x_1x_3)-3(x_1^2-x_2x_4),\\
Y_{33}&=& (x_2^2-3x_1x_3)+(x_1^2-x_2x_4),\\
Y_{23}&=& -x_1x_2+3x_3x_4.
\end{eqnarray*}
In the new coordinates on $T^*M^{(r)}\cong\C^4$, 
$L^{(r)}$ is the cone over the curve in $\P^3$ which is
defined in homogeneous coordinates by the vanishing of
$x_1^2-x_2x_4$, $x_2^2-3x_1x_3$ and $x_1x_2-3x_3x_4$. 
On the affine chart of $\P^3$ with $x_4=1$, this is the
twisted cubic $(x_1\mapsto (x_1,x_1^2,\frac{1}{3}x_1^3))$,
which is smooth. 
On the affine chart of $\P^3$ with $x_3=1$, this is the
smooth curve 
$(x_2\mapsto(\frac{1}{3}x_2^2,x_2,\frac{1}{9}x_2^3))$
which is also a twisted cubic.
The curve has no points with $x_3=x_4=0$. Therefore it is 
smooth. \hfill$\Box$

\begin{remark}\label{t6.6}
The second author is grateful to Paul Seidel how showed him
in 2000 this curve in $\P^3$ and explained that it is not 
a complete intersection, that it is a Legendre curve with 
respect to the 1-form $x_1\ddd x_2-x_2\ddd x_1-\ddd x_3$ 
(in the affine chart with $x_4=1$) and that the cone over it in 
$\C^4$ is smooth outside 0 and is Lagrange. 
\end{remark}

\section{Partial classification of 3-dimensional 
generically semisimple $F$-manifolds}\label{c7}
\setcounter{equation}{0}

The long Theorem \ref{t7.1} is the main result of
this section. It gives normal forms for all germs of
generically semisimple $F$-manifolds with
$T_{0}M\cong Q^{(2)}$, except $A_3$, $B_3$ and $H_3$. 
Part (a) of it is essentially Theorem 5.30
in \cite{He02}, but with some change in the normal form.
The parts (b)--(e) are new. 
Corollary \ref{t7.2} distinguishes those germs of $F$-manifolds
in Theorem \ref{t7.1} which have an Euler field.
The germs of $F$-manifolds in Theorem \ref{t7.1}
are closely related to the
germs of plane curves with multiplicity 3.
The Remarks \ref{t7.3} comment on this and make the 
cases in Theorem \ref{t7.1} more transparent.

\begin{theorem}\label{t7.1}
In the following, normal forms for all irreducible germs 
$(M,0)$of 3-dimensional generically semisimple $F$-manifolds
with $T_0M\cong Q^{(2)}$ except $A_3,B_3$ and $H_3$ 
are listed by their data in the Remarks \ref{t6.1}.  
Each isomorphism class of such a germ is represented 
by a finite positive number of normal forms. 
The normal forms split into 5 families with discrete
and holomorphic parameters, with

\begin{tabular}{c|c|c|c|c|c}
family in & (a) & (b) & (c) & (d) & (e) \\ \hline
number of components of $(L_M,\lambda)$ & 3 & 2 & 2 & 1 & 1\\ 
\hline
discrete parameters & $p,q$ & $p$ & $p,q$ & $p$ & $p$\\
\end{tabular}

\noindent
with $p,q\in\Z_{\geq 2}$ and $q\geq p$. 
There are always $p-1$ holomorphic parameters
$(\gamma_0,...,\gamma_{p-2})\in\C^{p-1}$ or in an open
subset. 
We use the notations in the Remarks \ref{t6.1},
especially $(M,0)=(\C^3,0)$ with coordinates $t=(t_1,t_2,t_3)$.
In all cases the caustic is $\KK=\{t\in M\,|\, t_2=0\}$.
It coincides with the $\mu$-constant stratum.
For $t\in\KK$ $T_tM\cong Q^{(2)}$. Locally on $M-\KK$,
the analytic spectrum is 
\begin{eqnarray}\label{7.1}
L_M=\bigcup_{j=1}^3\{(y,t)\in T^*M\,|\, 
y_1=1,y_2=\paa_2f^{(j)},y_3=h_2y_2^2+h_1y_2+h_0\}
\end{eqnarray}
with $h_2,h_1,h_0\in\C\{t_2,t_3\}$ as below,
with $h_2(0)\neq 0,h_1(0)=h_0(0)=0$. 
The Euler field on $M-\KK$
is $E=(t_1+c_1)\paa_1+\varepsilon_2\paa_2+\varepsilon_3\paa_3$
with $c_1\in \C$ and $\varepsilon_2,\varepsilon_3$ 
as below. Most often, $\varepsilon_3$ and $E$ 
are meromorphic along $\KK$ (see Corollary \ref{t7.2}
for the cases when they are holomorphic on $M$). 
The function
\begin{eqnarray}\label{7.2}
\rho:=t_2^{p-2}t_3+ \sum_{i=0}^{p-2}\gamma_it_2^i
\in \C\{t_2,t_3\}
\end{eqnarray}
will always turn up in some $f^{(j)}$. 

\medskip
(a) $(\gamma_0,...,\gamma_{p-2})\in\C^*\times \C^{p-2}$
with $\gamma_0\neq 1$ if $p=q$, \\
then $\rho\in \C\{t_2,t_3\}^*$, i.e. it is a unit in 
$\C\{t_2,t_3\}$, because $\gamma_0\neq 0$,\\
$f^{(1)}=0$, $f^{(2)}=t_2^p$, 
$f^{(3)}=t_2^q\cdot\rho$, \\
$L_M=\bigcup_{j=1}^3L^{(j)}$ has 3 smooth components,\\
$h_2=\Bigl((q+t_2\paa_2)(\rho)\Bigr)^{-1}\Bigl(((q+t_2\paa_2)(\rho))t_2^{q-p}
-p\Bigr)^{-1}\in\C\{t_2,t_3\}^*$, \\
$h_2^{-1}h_1=-pt_2^{p-1}$, $h_0=0$,\\
Euler field: $\varepsilon_2=\frac{1}{p}t_2$,  
$\varepsilon_3=-\frac{1}{p}t_2^{2-p}
((q-p+t_2\paa_2)(\rho))$.

\medskip
(b) $(\gamma_0,...,\gamma_{p-2})\in\C^{p-1}$, 
$\rho\in\C\{t_2,t_3\}$, \\
$f^{(1)}=0$, 
$f^{(2\&3)}=t_2^{\frac{1}{2}+p}+ t_2^{1+p}\cdot\rho$, \\
$L_M=L^{(1)}\cup L^{(2\&3)}$ has 1 smooth component
$L^{(1)}$ and 1 singular component $L^{(2\&3)}$,\\
$h_2=\Bigl((\frac{1}{2}+p)^{2}-t_2((1+p+t_2\paa_2)
(\rho))^2\Bigr)^{-1}
\in\C\{t_2,t_3\}^*$, \\
$h_2^{-1}h_1=-2 t_2^{p}((1+p+t_2\paa_2)(\rho))$, 
$h_0=0$,\\
Euler field: $\varepsilon_2=\frac{1}{\frac{1}{2}+p}t_2$,  
$\varepsilon_3=-\frac{1}{\frac{1}{2}+p}t_2^{2-p}
((\frac{1}{2}+t_2\paa_2)(\rho))$.

\medskip
(c) 
$(\gamma_0,...,\gamma_{p-2})\in\C^*\times\C^{p-2}$, 
and thus $\rho \in\C\{t_2,t_3\}^*$,\\
$f^{(1)}=0$, 
$f^{(2\&3)}=t_2^{\frac{1}{2}+q}\cdot\rho + t_2^{p}$,\\
$L_M=L^{(1)}\cup L^{(2\&3)}$ has 1 smooth component
$L^{(1)}$ and 1 singular component $L^{(2\&3)}$,\\
$h_2=\Bigl((\frac{1}{2}+q+t_2\paa_2)(\rho)\Bigr)^{-1} 
\Bigl(p-\frac{1}{p}t_2^{1+2(q-p)}
((\frac{1}{2}+q+t_2\paa_2)(\rho))^{2}\Bigr)^{-1}
\in\C\{t_2,t_3\}^*$, \\
$h_2^{-1}h_1=-pt_2^{p-1}-\frac{1}{p}t_2^{2q-p}
((\frac{1}{2}+q+t_2\paa_2)(\rho))^2$, $h_0=0$,\\
Euler field: $\varepsilon_2=\frac{1}{p}t_2$,  
$\varepsilon_3=-\frac{1}{p}t_2^{2-p}
((\frac{1}{2}+q-p+t_2\paa_2)(\rho))$.

\medskip
(d)  
$(\gamma_0,...,\gamma_{p-2})\in\C^{p-1}$, 
$\rho\in\C\{t_2,t_3\}$,\\
$f=f^{(1\&2\&3)}=t_2^{\frac{1}{3}+p}+ 
t_2^{\frac{2}{3}+p}\cdot\rho$,\\
$L_M$ is irreducible, \\
$h_2=\Bigl((\frac{1}{3}+p)^2 - 
t_2(\frac{1}{3}+p)^{-1} 
((\frac{2}{3}+p+t_2\paa_2)(\rho))^{3}\Bigr)^{-1}
\in\C\{t_2,t_3\}^*$, \\
$h_2^{-1}h_1=-t_2^{p}(\frac{1}{3}+p)^{-1}
((\frac{2}{3}+p+t_2\paa_2)(\rho))^2$,\\
$h_2^{-1}h_0=-2t_2^{2p-1}(\frac{1}{3}+p)
((\frac{2}{3}+p+t_2\paa_2)(\rho))$,\\
Euler field: $\varepsilon_2=\frac{1}{\frac{1}{3}+p}t_2$,  
$\varepsilon_3=-\frac{1}{\frac{1}{3}+p}t_2^{2-p}
((\frac{1}{3}+t_2\paa_2)(\rho))$.

\medskip
(e) 
$(\gamma_0,...,\gamma_{p-2})\in\C^{p-1}$, 
$\rho\in \C\{t_2,t_3\}$,\\
$f=f^{(1\&2\&3)}=t_2^{\frac{4}{3}+p}\cdot\rho +
t_2^{\frac{2}{3}+p}$,\\
$L_M$ is irreducible,\\
$h_2=\Bigl((\frac{2}{3}+p)^2 - 
t_2^2(\frac{2}{3}+p)^{-1} 
((\frac{4}{3}+p+t_2\paa_2)(\rho))^{3}\Bigr)^{-1}
\in\C\{t_2,t_3\}^*$, \\
$h_2^{-1}h_1=-t_2^{p+1}(\frac{2}{3}+p)^{-1}
((\frac{4}{3}+p+t_2\paa_2)(\rho))^2$,\\
$h_2^{-1}h_0=-2t_2^{2p}(\frac{2}{3}+p)
((\frac{4}{3}+p+t_2\paa_2)(\rho))$,\\
Euler field: $\varepsilon_2=\frac{1}{\frac{2}{3}+p}t_2$,  
$\varepsilon_3=-\frac{1}{\frac{2}{3}+p}t_2^{2-p}
((\frac{2}{3}+t_2\paa_2)(\rho))$.
\end{theorem}

{\bf Proof:}
The proofs of the parts (a)--(e) are similar.
Part (a) is essentially proved in \cite{He02}, though here
we chose a different normal form than in \cite{He02}. 
The first steps in the following proof hold for (a)--(e). 
We give all the details for part (b) and part (d).
We discuss differences and similarities for the parts
(c), (e) and (a). 

We consider an irreducible germ $(M,0)$ of a 
3-dimensional generically semisimple $F$-manifold
with $T_0M\cong Q^{(2)}$, which is not $A_3$, $B_3$ or $H_3$.
Theorem \ref{t6.3} says that then $(L_M,\lambda)\cong 
(\C^2,0)\times (C,0)$ where $(C,0)\subset(\C^2,0)$
is the germ of a plane curve with multiplicity 3.
And the caustic $\KK$ is isomorphic to the image in $M$
of the part of $L_M$ which is isomorphic to 
$(\C^2,0)\times\{0\}$, and $\KK$ is a smooth surface in $M$.

The coordinates $t=(t_1,t_2,t_3)$ can and will be chosen 
such that $\KK=\{t\in M\,|\, t_2=0\}$. 

$(C,0)$ has multiplicity 3, and therefore it has either 3
smooth components or 1 smooth and 1 singular component
or only 1 singular component. The corresponding
components of $L_M$ are called 
$L^{(j)}$, $j\in\{1,2,3\}$, in the first case,
$L^{(1)}$ and $L^{(2\&3)}$ in the second case 
and $L_M=L^{(1\&2\&3)}$ in the third case. 
The parts of the multivalued function $f$ 
which correspond to these components are called accordingly
$f^{(j)}$, $f^{(2\&3)}$ or $f^{(1\&2\&3)}$. 

Recall $F=t_1+f$ from the Notations \ref{t6.1}. 
The coordinate $t_1$ can and will be chosen 
(by a coordinate change as in \eqref{4.9}) such that 
\begin{eqnarray}\label{7.3}
f^{(1)}=0\quad\textup{ in the cases with 3 or 2 components,}\\
f^{(1)}+f^{(2)}+f^{(3)}=0\quad
\textup{ in the cases with 1 component.}\label{7.4}
\end{eqnarray}
In fact, in all cases $f^{(1)}+f^{(2)}+f^{(3)}$ is univalued,
and $t_1$ can be chosen such that \eqref{7.4} holds.
Then also $\paa_2(f^{(1)}+f^{(2)}+f^{(3)})=0$
and $\paa_3(f^{(1)}+f^{(2)}+f^{(3)})=0$. 
We see that this choice was already discussed in 
Remark \ref{t4.5}.
In the generically semisimple case, the function $F$
gives an alternative starting point for understanding
this choice of the coordinate $t_1$. 

In the cases with 1 component, we use this choice in 
\eqref{7.4}, and there it gives
\begin{eqnarray}
\prod_{j=1}^3(x-\paa_2f^{(j)}) =x^3+g_1x+g_0,
\quad \textup{ so }g_2=0.\label{7.5}
\end{eqnarray}

In the cases with 3 or 2 components, we prefer 
the choice in \eqref{7.3}, as it makes there the 
calculations easier. 
Then in the cases with 3 or 2 components, \eqref{6.6}
for $j=1$ gives $h_0=0$ and 
\begin{eqnarray}
\prod_{j=1}^3(x-\paa_2f^{(j)}) =x^3+g_2x^2+g_1x
\quad\textup{with}
\nonumber\\
g_2=-\paa_2f^{(2)}-\paa_2f^{(3)},\ 
g_1=\paa_2f^{(2)}\cdot\paa_2f^{(3)}.\label{7.6}
\end{eqnarray}

(b) and (c) 
Now we turn to the cases where $L_M$ has 2 components,
the smooth component $L^{(1)}$ and the singular
component $L^{(2\&3)}$. We have $f^{(1)}=0$ and 
\begin{eqnarray}\label{7.7}
f^{(2\&3)}=t_2^{1/2+p_1}\rho_1 + t_2^{p_2}\rho_2
\end{eqnarray}
with $\rho_1\in\C\{t_2,t_3\}-t_2\C\{t_2,t_3\}$ and
$\rho_2\in\C\{t_2,t_3\}-(t_2\C\{t_2,t_3\}-\{0\})$.
Here $\rho_1,\rho_2$ and $p_1\in\Z_{\geq 0}$ are unique,
and $p_2\in\Z_{\geq 0}$ is unique if $\rho_2\neq 0$.
If $\rho_2=0$, we put $p_2:=\infty$. 

The branched covering $\pi_L:L_M\to M$ is branched only over
$\KK=\{t\in M\,|\, t_2=0\}$. This implies two facts: 
First, $L^{(1)}$ and $L^{(2\&3)}$ intersect only over
$\KK$, and second, the branched covering $\pi_L:L^{(2\&3)}\to M$
is branched only over $\KK$. 
The second fact tells $\rho_1\in\C\{t_2,t_3\}^*$,
i.e. $\rho_1$ is a unit, i.e. $\rho_1(0)\neq 0$. 
If $p_1< p_2$, this is sufficient also for the first fact.
Then we are in the cases in (b).
If $p_1\geq p_2$, the first fact tells 
$\rho_2\in\C\{t_2,t_3\}^*$.
Then we are in the cases in (c). 

\medskip
(b) Now we turn to the cases in (b), i.e. $f^{(1)}=0$ and 
$f^{(2\&3)}$ as in \eqref{7.7} with $p_1<p_2$. 
Rename $p:=p_1$. Then $t_2$ can and will be chosen such 
that $t_2^{\frac{1}{2}+p}\rho_1=t_2^{\frac{1}{2}+p}$.
Then we write 
\begin{eqnarray}\label{7.8}
f^{(2\&3)}=t_2^{\frac{1}{2}+p}+t_2^{1+p} \rho
\end{eqnarray}
for some $\rho\in\C\{t_2,t_3\}$. Next we will exploit
\eqref{6.6} together with $h_2(0)\neq 0$ in order to put $\rho$
into a normal form by a good choice of $t_3$, and to 
calculate $h_2$ and $h_1$ (recall $h_0=0$ because of 
\eqref{6.6} for $f^{(1)}=0$). \eqref{6.6} gives
\begin{eqnarray}\label{7.9}
t_2^{1+p}\paa_3\rho &=&\paa_3f^{(2\&3)}= 
h_2\cdot \paa_2 f^{(2\&3)}
\cdot \Bigl(\paa_2f^{(2\&3)}+h_2^{-1}h_1\Bigr)\\
&=& h_2\Bigl(\bigl[(\paa_2t_2^{\frac{1}{2}+p})^2+
\paa_2(t_2^{1+p}\rho)(\paa_2(t_2^{1+p}\rho)
+h_2^{-1}h_1)\bigr]\nonumber\\
&&+\bigl[\paa_2t_2^{\frac{1}{2}+p}
(2\paa_2(t_2^{1+p}\rho)+h_2^{-1}h_1)\bigr]\Bigr).\label{7.10}
\end{eqnarray}
The term in square brackets in the line \eqref{7.10}
must vanish because of the half-integer exponent of $t_2$.
This allows to calculate
$h_2^{-1}h_1=-2\paa_2(t_2^{1+p}\rho)$, 
see the formula in part (b) in the theorem.
And it simplifies the other summand, 
\begin{eqnarray}
\paa_3\rho&=& 
h_2t_2^{-1-p}\Bigl((\paa_2t_2^{\frac{1}{2}+p})^2-
(\paa_2(t_2^{1+p}\rho))^2\Bigr)\nonumber\\
&=&h_2 t_2^{p-2}\Bigl((\frac{1}{2}+p)^2
-t_2((1+p+t_2\paa_2)(\rho))^2\Bigr),\label{7.11}
\end{eqnarray}
so $\paa_3\rho
=t_2^{p-2}\cdot(\textup{a unit in }\C\{t_2,t_3\})$.
This implies $p\in\Z_{\geq 2}$. And we can and will choose
$t_3$ such that $\rho$ is as in \eqref{7.2}. Then 
$h_2$ is determined by \eqref{7.11} with $\paa_3\rho=t_2^{p-2}$. 

The coefficients $\varepsilon_2$ and $\varepsilon_3$ of the
Euler field are determined by \eqref{6.3} for $f^{(2\&3)}$
as in \eqref{7.8} and \eqref{7.2}: 
\begin{eqnarray}\label{7.12}
&&t_2^{\frac{1}{2}+p}+t_2^{1+p}\rho = f^{(2\&3)}
= \varepsilon_2\paa_2 f^{(2\&3)}+\varepsilon_3\paa_3f^{(2\&3)}\\
&=&\varepsilon_2\Bigl((\frac{1}{2}+p)t_2^{-\frac{1}{2}+p}
+((1+p+t_2\paa_2)(\rho))t_2^{p}\Bigr)
+\varepsilon_3 t_2^{2p-1}.\nonumber
\end{eqnarray}
Comparison of the terms with half-integer exponents gives
$\varepsilon_2=(\frac{1}{2}+p)^{-1}t_2$. Then comparison
of the terms with integer exponents gives $\varepsilon_3$.

Almost all steps in this reduction process to a normal
form were unique. \eqref{7.7} was the general ansatz.
The choice of $t_2$ with $t_2^{\frac{1}{2}+p}\rho_1
=t_2^{\frac{1}{2}+p}$ was unique up to a unit root of order
$1+2p$. The choice of $t_3$ was unique. 
Therefore the isomorphism class of $(M,0)$ is represented
by up to $1+2p$ normal forms.

\medskip
(c) Now we turn to the cases in (c), i.e. $f^{(1)}=0$ and 
$f^{(2\&3)}$ as in \eqref{7.7} with $p_1\geq p_2$. 
Rename $q:=p_1$ and $p:=p_2$. 
Above we showed that $\rho_1$ and $\rho_2$ are units in
$\C\{t_2,t_3\}$. 
We can and will choose $t_2$ such that 
$t_2^{p}\rho_2=t_2^{p}$. Then we write 
\begin{eqnarray}\label{7.13}
f^{(2\&3)}=t_2^{\frac{1}{2}+q} \rho+t_2^{p}
\end{eqnarray}
for some $\rho\in\C\{t_2,t_3\}^*$.
As in the proof of part (b), the next step is to exploit
\eqref{6.6} together with $h_2(0)\neq 0$ in order to
put $\rho$ into a normal form by a good choice of $t_3$,
and to calculate $h_2$ and $h_1$. The calculation is similar
to the calculation of \eqref{7.9} above. It leads to 
$\paa_3\rho=t_2^{p-2}\cdot(\textup{a unit in }\C\{t_2,t_3\})$.
This implies $p\in\Z_{\geq 2}$. And it allows to choose
$t_3$ such that $\rho$ is as in \eqref{7.2}. We skip the details
of the calculations. The results are written in part (c) 
in the theorem.
The fact that here $\rho$ is a unit, implies 
$\gamma_0\in\C^*$. Also the calculation of the coefficients
$\varepsilon_2$ and $\varepsilon_3$ of the Euler field
is similar to the calculation \eqref{7.12} above.
Again we skip the details. The results are written in part (c)
in the theorem. The choice of $t_2$ was unique up to 
a unit root of order $p$. The choice of $t_3$ was unique.
Therefore the isomorphism class of $(M,0)$ is represented
by up to $p$ normal forms. 

\medskip
(d) and (e) Now we turn to the cases where $L_M$ is irreducible.
A priori we have
\begin{eqnarray}\label{7.14}
f=t_2^{\frac{1}{3}+p_1}\rho_1+t_2^{\frac{2}{3}+p_2}\rho_2
+t_2^{1+p_3}\rho_3
\end{eqnarray}
with $\rho_1,\rho_2,\rho_3\in\C\{t_2,t_3\}
-(t_2\C\{t_2,t_3\}-\{0\})$ and $(\rho_1,\rho_2)\neq (0,0)$. 
But $\sum_{j=1}^3f^{(j)}=0$ tells $\rho_3=0$. 
If $\rho_1\neq 0$ then $\rho_1$ and $p_1\in\Z_{\geq 0}$
are unique, else $p_1:=\infty$. 
If $\rho_2\neq 0$ then $\rho_2$ and $p_2\in\Z_{\geq 0}$
are unique, else $p_2:=\infty$. 

The branched covering $\pi_L:L_M\to M$ is branched
only over $\KK=\{t\in M\,|\, t_2=0\}$. If $p_1\leq p_2$,
this implies $\rho_1\in \C\{t_2,t_3\}^*$,
and then we are in the cases in (d). 
If $p_1> p_2$,
this implies $\rho_2\in \C\{t_2,t_3\}^*$,
and then we are in the cases in (e). 

\medskip
(d) Now we turn to the cases in (d), i.e. $f$ is as in 
\eqref{7.14} with $\rho_3=0$ and 
$p_1\leq p_2$ and $\rho_1\in\C\{t_2,t_3\}^*$.
Rename $p:=p_1$. 
Then $t_2$ can and will be chosen such that 
$t_2^{\frac{1}{3}+p}\rho_1=t_2^{\frac{1}{3}+p}$.
Then we write 
\begin{eqnarray}\label{7.15}
f=t_2^{\frac{1}{3}+p}+t_2^{\frac{2}{3}+p} \rho
\end{eqnarray}
for some $\rho\in\C\{t_2,t_3\}$.
As in the proofs of the parts (b) and (c), the next step is
to exploit \eqref{6.6} together with $h_2(0)\neq 0$ in order
to put $\rho$ into a normal form by a good choice of $t_3$,
and to calculate $h_2$, $h_1$ and $h_0$. The calculation is
as follows.
\begin{eqnarray}\label{7.16}
t_2^{\frac{2}{3}+p}\paa_3\rho 
&=& \paa_3f = h_2\cdot \Bigl((\paa_2 f)^2
+h_2^{-1}h_1\paa_2f+h_2^{-1}h_0\Bigr)\\
&=& h_2\Bigl(\bigl[(\paa_2t_2^{\frac{1}{3}+p})^2 
+ h_2^{-1}h_1\paa_2(t_2^{\frac{2}{3}+p}\rho)\bigr]\label{7.17}\\
&+& \bigl[(\paa_2(t_2^{\frac{2}{3}+p}\rho))^2
+h_2^{-1}h_1\paa_2t_2^{\frac{1}{3}+p}\bigr]\label{7.18}\\
&+& \bigl[2\paa_2t_2^{\frac{1}{3}+p}
\paa_2(t_2^{\frac{2}{3}+p}\rho) +h_2^{-1}h_0\bigr]\Bigr)
\label{7.19}
\end{eqnarray}
The terms in square brackets in the lines \eqref{7.18}
and \eqref{7.19} must vanish because of the 
exponents in $\frac{1}{3}+\Z$ and $\Z$ of $t_2$.
This allows to calculate $h_2^{-1}h_0$ and 
$h_2^{-1}h_1=-(\paa_2t_2^{\frac{1}{3}+p})^{-1}
(\paa_2(t_2^{\frac{2}{3}+p}\rho))^2$, 
see the formulas in part (d) in the theorem.
And it simplifies the term in square brackets in the line
\eqref{7.17}, 
\begin{eqnarray}
\paa_3\rho&=& 
h_2t_2^{-\frac{2}{3}-p}\Bigl((\paa_2t_2^{\frac{1}{3}+p})^2-
(\paa_2t_2^{\frac{1}{3}+p})^{-1}
(\paa_2(t_2^{\frac{2}{3}+p}\rho))^3\Bigr)\nonumber\\
&=&h_2 t_2^{p-2}\Bigl((\frac{1}{3}+p)^2
-t_2(\frac{1}{3}+p)^{-1}((\frac{2}{3}+p+t_2\paa_2)(\rho))^3\Bigr),\label{7.20}
\end{eqnarray}
so $\paa_3\rho
=t_2^{p-2}\cdot(\textup{a unit in }\C\{t_2,t_3\})$.
This implies $p\in\Z_{\geq 2}$. And we can and will choose
$t_3$ such that $\rho$ is as in \eqref{7.2}. Then $h_2$ 
is determined by \eqref{7.20} with $\paa_3\rho=t_2^{p-2}$. 

The coefficients $\varepsilon_2$ and $\varepsilon_3$ of the
Euler field are determined by \eqref{6.3} for $f$
as in \eqref{7.15} and \eqref{7.2}: 
\begin{eqnarray}\label{7.21}
&&t_2^{\frac{1}{3}+p}+t_2^{\frac{2}{3}+p}\rho = f
= \varepsilon_2\paa_2 f+\varepsilon_3\paa_3f\\
&=&\varepsilon_2\Bigl((\frac{1}{3}+p)t_2^{-\frac{2}{3}+p}
+((\frac{2}{3}+p+t_2\paa_2)(\rho))t_2^{-\frac{1}{3}+p}\Bigr)
+\varepsilon_3 t_2^{-\frac{4}{3}+2p}.\nonumber
\end{eqnarray}
Comparison of the terms with exponents in $\frac{1}{3}+\Z$ gives
$\varepsilon_2=(\frac{1}{3}+p)^{-1}t_2$. Then comparison
of the terms with exponents in $\frac{2}{3}+\Z$ 
gives $\varepsilon_3$.
The choice of $t_2$ with $t_2^{\frac{1}{3}+p}\rho_1
=t_2^{\frac{1}{3}+p}$ was unique up to a unit root of order
$1+3p$. The choice of $t_3$ was unique. 
Therefore the isomorphism class of $(M,0)$ is represented
by up to $1+3p$ normal forms.

\medskip
(e) Now we turn to the cases in (e), i.e. $f$ as in \eqref{7.15}
with $p_1>p_2$ and $\rho_2\in\C\{t_2,t_3\}^*$. 
Rename $p:=p_2$. Then $t_2$ can and will be chosen such that 
$t_2^{\frac{2}{3}+p}\rho_2=t_2^{\frac{2}{3}+p}$. Then we write
\begin{eqnarray}\label{7.22}
f=t_2^{\frac{4}{3}+p} \rho+t_2^{\frac{2}{3}+p}
\end{eqnarray}
for some $\rho\in\C\{t_2,t_3\}$.
As in the proofs of the parts (b), (c) and (d), 
the next step is to exploit
\eqref{6.6} together with $h_2(0)\neq 0$ in order to
put $\rho$ into a normal form by a good choice of $t_3$,
and to calculate $h_2$, $h_1$ and $h_0$. 
The calculation is similar
to the calculation of \eqref{7.16} above. It leads to 
$\paa_3\rho=t_2^{p-2}\cdot(\textup{a unit in }\C\{t_2,t_3\})$.
This implies $p\in\Z_{\geq 2}$. And it allows to choose
$t_3$ such that $\rho$ is as in \eqref{7.2}. We skip the details
of the calculations. The results are written in part (e) 
in the theorem.
Also the calculation of the coefficients
$\varepsilon_2$ and $\varepsilon_3$ of the Euler field
is similar to the calculation \eqref{7.21} above.
Again we skip the details. The results are written in part (e)
in the theorem. The choice of $t_2$ was unique up to 
a unit root of order $2+3p$. The choice of $t_3$ was unique.
Therefore the isomorphism class of $(M,0)$ is represented
by up to $2+3p$ normal forms. 

\medskip
(a) Now we turn to the cases in (a), the cases where
$L$ has 3 components. They were treated in Theorem 5.30
in \cite{He02}. But here we choose the normal forms
a bit differently. 

The plane curve germs
\begin{eqnarray}\label{7.23}
(C^{(j)},0):=\{(y_2,t_2)\in (\C^2,0)
\,|\, y_2=\paa_2f^{(j)}(t_2,0)\}
\end{eqnarray}
in the $(y_2,t_2)$-plane satisfy 
$(L,\lambda)\cong (\C^2,0)\times \bigcup_{j=1}^3(C^{(j)},0)$ 
by the proof of Theorem \ref{t6.3} (a). 
We choose their numbering such
that the pair $(C^{(1)},C^{(3)})$ has the highest intersection 
number, which we call $q-1$ for some $q\in\Z_{\geq 2}$.
Then the pairs $(C^{(1)},C^{(2)})$ and $(C^{(2)},C^{(3)})$
have the same intersection number $p-1$ for some 
$p\in\Z_{\geq 2}$ with $p\leq q$. 

We have $f^{(1)}=0$ by \eqref{7.3} and 
\begin{eqnarray}\label{7.24}
f^{(2)}=t_2^{p_1}\rho_1,\quad f^{(3)}=t_2^{p_2}\rho_2
\end{eqnarray}
with $\rho_1,\rho_2\in\C\{t_2,t_3\}-t_2\C\{t_2,t_3\}$
and $p_1,p_2\in\N$.

The branched covering $\pi_L:L_M\to M$ is branched only over
$\KK=\{t\in M\,|\, t_2=0\}$, so the components 
$L^{(i)}$ and $L^{(j)}$ of $L_M$ intersect only over $\KK$. 
This and $(L,\lambda)\cong (\C^2,0)\times \bigcup_{j=1}^3(C^{(j)},0)$ 
shows $\rho_1,\rho_2\in\C\{t_2,t_3\}^*$, $p_1=p$, $p_2=q$, 
and in the case $p=q$ additionally $\rho_1(0)\neq \rho_2(0)$. 

$t_2$ can and will be chosen
such that $t_2^p\rho_1=t_2^p$. Then we have
\begin{eqnarray}\label{7.25}
f^{(1)}=0,\quad f^{(2)}=t_2^p,\quad f^{(3)}=t_2^q\cdot\rho
\end{eqnarray}
for some $\rho\in\C\{t_2,t_3\}^*$ with $\rho(0)\neq 1$ 
if $p=q$. As in the proofs of the parts (b)--(e), the next
step is to exploit \eqref{6.6} together with $h_2(0)\neq 0$
in order to put $\rho$ into a normal form by a good
choice of $t_3$, and to calculate $h_2$ and $h_1$.
The calculation is similar to the calculations in (b)--(e),
and, in fact, easier. It allows to choose $t_3$ such that
$\rho$ is as in \eqref{7.2}.
Then $h_2$ and $h_1$ are as in the theorem. Also the 
calculation of the Euler field is similar to the calculations
in (b)--(e). The numbering of $L^{(1)},L^{(2)}$ and $L^{(3)}$
was unique up to a permutation of $L^{(1)}$ and $L^{(3)}$
if $p<q$ and arbitrary if $p=q$. The choice of $t_2$
such that $t_2^p\rho_1=t_2^p$ was unique up to a unit
root of order $p$. The choice of $t_3$ was unique. 
Therefore the isomorphism class of $(M,0)$ is represented 
by up to $2p$ or $6p$ normal forms.
\hfill$\Box$

\bigskip

Probably the most interesting of the $F$-manifolds in
Theorem \ref{t7.1} are those where the Euler field
is holomorphic on $M$. The next corollary makes them explicit. 

\begin{corollary}\label{t7.2}
Each irreducible germ $(\www M,0)$ of a 3-dimensional 
generically semisimple $F$-manifold with 
$T_0\www M\cong Q^{(2)}$ and with (holomorphic) Euler field 
is isomorphic to one of the germs $A_3$, $B_3$ or $H_3$ or to
a germ $(M,(t_1,0,t_3))$ for suitable $t_1\in\C$ and $t_3\in\C$ 
(or $\C^*$ or $\C-\{0;1\}$) of one of the
$F$-manifolds which are listed below. We use the same
notations as in Theorem \ref{t7.1}. There are 7 families
of $F$-manifolds. The family in (a)(i) has no parameter, so 
there is a single $F$-manifold. 
The family in (a)(iii) has one discrete and one holomorphic
parameter $\gamma_0$. 
The other families have one discrete parameter
and no holomorphic parameter.The Euler field is 
$E=t_1\paa_1+\varepsilon_2\paa_2+\varepsilon_3\paa_3$
with $\varepsilon_2,\varepsilon_3$ as below. 

\noindent
\begin{tabular}{c|c|c|c|c|c|c|c}
family in & (a)(i) & (a)(ii) & (a)(iii) & (b) & (c) & (d) 
& (e) \\ \hline
number of comp. of $(L_M,\lambda)$ 
& 3 & 3 & 3 & 2 & 2 & 1 & 1\\ 
\hline
discrete parameter & -- & $q$ & $p$ & $p$ & $q$ & $p$ & $p$\\
\end{tabular}

\medskip
(a) (i) $M=\C^2\times(\C-\{0;1\})$, ($p=q=2$,)\\
$f^{(1)}=0$, $f^{(2)}=t_2^2$, $f^{(3)}=t_2^2t_3$, \\
$h_2=(4t_3(t_3-1))^{-1}$, 
$h_2^{-1}h_1=-2t_2$, $h_0=0$,\\
Euler field: $\varepsilon_2=\frac{1}{2}t_2$, 
$\varepsilon_3=0$. 

\medskip
(a) (ii) $M=\C^2\times\C^*$, ($p=2$,) $q\in\Z_{\geq 3}$,\\
$f^{(1)}=0$, $f^{(2)}=t_2^2$, 
$f^{(3)}=t_2^qt_3$, \\
$h_2=(qt_3)^{-1}(qt_3t_2^{q-2}-2)^{-1}$, 
$h_2^{-1}h_1=-2t_2$, $h_0=0$,\\
Euler field: $\varepsilon_2=\frac{1}{2}t_2$, 
$\varepsilon_3=-\frac{q-2}{2}t_3$. 

\medskip
(a) (iii) $M=\C^3$, ($p=q$,) $p\in\Z_{\geq 3}$,
$\gamma_0\in\C-\{0;1\}$, \\
$f^{(1)}=0$, $f^{(2)}=t_2^p$, 
$f^{(3)}=t_2^p(\gamma_0+t_2^{p-2}t_3)$, \\
$h_2=(p\gamma_0+(p-2)t_2^{p-2}t_3)^{-1}
(p(\gamma_0-1)+(p-2)t_2^{p-2}t_3)^{-1}$, \\
$h_2^{-1}h_1=-pt_2^{p-1}$, $h_0=0$,\\
Euler field: $\varepsilon_2=\frac{1}{p}t_2$, 
$\varepsilon_3=-\frac{p-2}{p}t_3$. 

\medskip
(b) $M=\C^3$, $p\in\Z_{\geq 2}$,\\
$f^{(1)}=0$, 
$f^{(2\&3)}=t_2^{\frac{1}{2}+p}+ t_2^{2p-1}t_3$, \\
$h_2=((\frac{1}{2}+p)^{2}-(2p-1)^2t_2^{2p-3}t_3^2)^{-1}$, \\
$h_2^{-1}h_1=-2(2p-1)t_2^{2p-2}t_3$, 
$h_0=0$,\\
Euler field: $\varepsilon_2=\frac{1}{\frac{1}{2}+p}t_2$,  
$\varepsilon_3=-\frac{p-\frac{3}{2}}{\frac{1}{2}+p}t_3$.

\medskip
(c) $M=\C^2\times\C^*$, ($p=2$,) $q\in\Z_{\geq 2}$,\\
$f^{(1)}=0$, 
$f^{(2\&3)}=t_2^{\frac{1}{2}+q}t_3 + t_2^{2}$,\\
$h_2=((\frac{1}{2}+q)t_3)^{-1} 
\Bigl(2-\frac{1}{2}(\frac{1}{2}+q)^2t_2^{2q-3}t_3^2
\Bigr)^{-1}$, \\
$h_2^{-1}h_1=-2t_2-\frac{1}{2}(\frac{1}{2}+q)^2t_2^{2q-2}
t_3^2$, $h_0=0$,\\
Euler field: $\varepsilon_2=\frac{1}{2}t_2$,  
$\varepsilon_3=-\frac{1}{2}
(q-\frac{3}{2})t_3$.

\medskip
(d)  $M=\C^3$, $p\in\Z_{\geq 2}$,\\
$f=f^{(1\&2\&3)}=t_2^{\frac{1}{3}+p}+ 
t_2^{2p-\frac{4}{3}}t_3$,\\
$h_2=\Bigl((\frac{1}{3}+p)^2 - (\frac{1}{3}+p)^{-1} 
((2p-\frac{4}{3})^3t_2^{3p-5}t_3^3\Bigr)^{-1}$, \\
$h_2^{-1}h_1=-(\frac{1}{3}+p)^{-1}
(2p-\frac{4}{3})^2t_2^{3p-4}t_3^2$,\\
$h_2^{-1}h_0=-2(\frac{1}{3}+p)
(2p-\frac{4}{3})t_2^{3p-3}t_3$,\\
Euler field: $\varepsilon_2=\frac{1}{\frac{1}{3}+p}t_2$,  
$\varepsilon_3=-\frac{p-\frac{5}{3}}{\frac{1}{3}+p}t_3.$

\medskip
(e) $M=\C^3$, $p\in\Z_{\geq 2}$, \\
$f=f^{(1\&2\&3)}=t_2^{2p-\frac{2}{3}}t_3 +
t_2^{\frac{2}{3}+p}$,\\
$h_2=\Bigl((\frac{2}{3}+p)^2 - 
(\frac{2}{3}+p)^{-1} 
(2p-\frac{2}{3})^3t_2^{3p-4}t_3^3\Bigr)^{-1}$, \\
$h_2^{-1}h_1=-(\frac{2}{3}+p)^{-1}
(2p-\frac{2}{3})^2t_2^{3p-3}t_3^2$,\\
$h_2^{-1}h_0=-2(\frac{2}{3}+p)
(2p-\frac{2}{3})t_2^{3p-2}t_3$,\\
Euler field: $\varepsilon_2=\frac{1}{\frac{2}{3}+p}t_2$,  
$\varepsilon_3=-\frac{p-\frac{4}{3}}{\frac{2}{3}+p}t_3$.
\end{corollary}

{\bf Proof:}
The shape of the Euler field in Theorem \ref{t6.3}
tells precisely under which conditions it is holomorphic.
The conditions are as follows.

(a) $p=2$ or ($p=q$ and $\gamma_1=...=\gamma_{p-3}=0$). 

(b) $\gamma_0=...=\gamma_{p-3}=0$.

(c) $p=2$. 

(d) $\gamma_0=...=\gamma_{p-3}=0$.

(e) $\gamma_0=...=\gamma_{p-3}=0$.

In all cases, we consider germs also at points with $t_3\neq0$, 
and therefore we can replace $\gamma_{p-2}+t_3$ by $t_3$ in 
$\rho$. A condition on $\gamma_{p-2}$
(to be in $\C^*$ or $\C-\{0;1\}$) translates into a condition
on $t_3$. This gives all statements in the corollary.
\hfill$\Box$ 

\begin{remarks}\label{t7.3} 
(i) The classification in Theorem \ref{t7.1} of 3-dimensional 
germs $(M,0)$ of generically semisimple $F$-manifolds with 
$T_0M\cong Q^{(2)}$ which are different from $A_3,B_3,H_3$ 
is precise, but not so transparent. 
It becomes more transparent if one takes a closer look at
the reduced plane curve germs $(C,0)$ with 
$(L_M,\lambda)\cong (\C^2,0)\times (C,0)$.
By Theorem \ref{t6.3}, they have multiplicity 3.
And by Corollary \ref{t4.7}, all reduced
plane curve germs $(C,0)$ with multiplicity 3 appear. 

\medskip
(ii) Each reduced plane curve germ has a topological
type. See \cite[3.4]{GLS07} for its definition. 
An old result of Brauner and Zariski (see e.g. 
\cite[Lemma 3.31 + Proposition 3.41 + Theorem 3.42]{GLS07}) 
is that the topological type of a reduced plane curve germ
is determined by the topological types of the irreducible
components and by their intersection numbers.
And the topological type of an irreducible plane curve germ
is determined by its Puiseux pairs 
(see e.g. \cite[3.4]{GLS07} for their definition).
The topological types of reduced plane curve germs of 
multiplicity 3 can be described and listed as follows.
In all cases, the number $p\in\Z_{\geq 2}$ and, if it
exists, also the number $q\in\Z_{\geq 2}$ are topological
invariants. 
\begin{list}{}{}
\item[(a)]
3 smooth curve germs $C^{(1)},C^{(2)},C^{(3)}$ with 
intersection multiplicities 
$i(C^{(1)},C^{(2)})=p-1$, $i(C^{(1)},C^{(3)})=q-1$, 
$i(C^{(2)},C^{(3)})=p-1$ for $p,q\in\Z_{\geq 2}$ with
$q\geq p$. 
\item[(b)]
1 smooth germ $C^{(1)}$ and one germ $C^{(2\&3)}$ of type 
$A_{2p-2}$ (namely with normal form $x_1^{2p-1}+x_2^2$) with the 
maximal possible intersection number 
$2p-1=i(C^{(1)},C^{2\&3)})$. 
\item[(c)]
1 smooth germ $C^{(1)}$ and one germ $C^{(2\&3)}$ of type 
$A_{2q-2}$ with an even intersection number 
$i(C^{(1)},C^{(2\&3)})=2p-2$ for $q,p\in\Z_{\geq 2}$
with $q\geq p$.  
\item[(d)]
1 irreducible germ with the only Puiseux pair
$(3p-2,3)$, so with a parametrization 
$(x\mapsto (x,\sum_{n\geq 3p-2}a_nx^{n/3}))$ with $a_n\in\C$
and $a_{3p-2}\neq 0$. 
\item[(e)]
1 irreducible germ with the only Puiseux pair
$(3p-1,3)$, so with a parametrization 
$(x\mapsto (x,\sum_{n\geq 3p-1}a_nx^{n/3}))$ with $a_n\in\C$
and $a_{3p-1}\neq 0$. 
\end{list}
One sees that the cases (a)--(e) correspond precisely
to the cases (a)--(e) in Theorem \ref{t7.1}. 
There the curve $(C,0)$ is the zero set of the 
polynomial $\prod_{j=1}^3(y_2-\paa_2f^{(j)}|_{t_3=0})
\in\C[y_2,t_2]$. 

\medskip
(iii) The following topological types contain quasihomogeneous
plane curve germs $(C,0)$: all in (b), (d) and (e);
those in (a) with $p=2$ or $p=q$; those in (c) with $p=2$. 
This fits to the cases in Corollary \ref{t7.2}.
The topological types in (a) with $p=q\geq 3$ contain
a 1-parameter family of quasihomogeneous curves (up to
coordinate changes). This gives the holomorphic parameter 
$\gamma_0$ in Corollary \ref{t7.2} (a) (iii). 

\medskip
(iv)  In all cases in Theorem \ref{t7.1},
there are $p-1$ holomorphic parameters 
$(\gamma_0,\gamma_1,...,\gamma_{p-3},\gamma_{p-2}+t_3)$
for the germs of $F$-manifolds. 
Here the last parameter $\gamma_{p-2}+t_3$ is an 
{\it internal parameter}, it is the parameter of the
1-dimensional $\mu$-constant stratum. 
The other parameters $(\gamma_1,...,\gamma_{p-3})$
(for $p\geq 3$; no other parameter for $p=2$) 
catch the isomorphism class of the plane curve germ 
$(C,0)$ and the choice of a symplectic structure
on the germ $(\C^2,0)$ of the $(y_2,t_2)$-plane.
This is the choice of a volume form, i.e. a form
$u\ddd y_2\ddd t_2$ with $u\in\C\{y_2,t_2\}^*$.
In \cite[Remark 5.31]{He02} the 3 types of parameters
are rephrased as follows.
\begin{list}{}{}
\item[($\alpha$)]
moduli for the complex structure of the germ $(C,0)$,
\item[($\beta$)]
moduli for the Lagrange structure of $(C,0)$ or, equivalently,
for the symplectic structure of $(\C^2,0)\supset (C,0)$,
\item[($\gamma$)]
moduli for the Lagrange fibration.
\end{list}
Here we have only one parameter of type $(\gamma)$, the 
internal parameter $\gamma_{p-2}+t_3$.
Remarkably, the sum of the numbers of parameters of types 
$(\alpha)$ and $(\beta)$ is constant, it is $p-2$. 
This is remarkable, as the number of parameters
of type $(\alpha)$ depends on the plane curve germ $(C,0)$
with which one starts. It has the shape $\tau(C,0)-b^{top}$,
where $b^{top}\in\N$ is a topological invariant and 
the {\it Tjurina number} $\tau(C,0)$ was defined in Theorem 
\ref{t2.16} (b) and is not a topological invariant. 
But by Theorem \ref{t2.16} (b)+(c), the number of parameters
of type ($\beta$) compensates this, as it is precisely
$\mu-\tau(C,0)=\dim H^1_{Giv}(\C^2,C,0)$. 
So the sum of the numbers of parameters of types
$(\alpha)$ and $(\beta)$ is $\mu-b^{top}$.
Here this number is $p-2$. 

\medskip
(v) The same number $\mu-b^{top}$ is also the number
of parameters of right equivalence classes of plane
curve germs with fixed topological type. Here the
right equivalence class is the class up to holomorphic
coordinate changes. 
This follows from the fact that $\mu-\tau(C,0)$ is also
the difference of the dimensions of the base space of a 
universal unfolding of a function germ for $(C,0)$ 
and of a semiuniversal deformation of $(C,0)$. 
But for a given reduced plane curve germ $(C,0)$, 
there is no canonical relation between the choice of a 
volume form on $(\C^2,0)$ and the choice of a function 
germ $f:(\C^2,0)\to(\C,0)$ with $(f^{-1}(0),0)=(C,0)$. 

\medskip
(vi) Here the normal forms in Theorem \ref{t7.1} are
misleading. In all topological types which contain
quasihomogeneous curves (up to coordinate changes),
the following holds (and probably it holds also for the
other topological types in (a) and (c)): 
The parameters $(\gamma_0,...,\gamma_{p-3})$ in 
$\prod_{j=1}^3(y_2-\paa_2f^{(j)}|_{t_3=0})$ are also
the parameters for the right equivalence classes.
And if they are fixed, the internal parameter
$\gamma_{p-2}+t_3$ does not change the right equivalence
class. This follows by inspection of the curves and
a description of the $\mu$-constant stratum in a 
universal unfolding of a quasihomogeneous singularity in
\cite{Va82}.

\medskip
(vii) The property in (vi), that the internal parameter
$\gamma_{p-2}+t_3$ does not change the right equivalence class, 
is a lucky coincidence of the chosen normal forms. 
It is easy to construct a concrete 
description of a germ of an $F$-manifold in Theorem 
\ref{t7.1} where this does not hold. 
Start with a plane curve germ $(C,0)$ which is not 
quasihomogeneous (up to coordinate changes) 
and choose function germs 
$g_2^{(0)}:=0$ and $g_1^{(0)},g_0^{(0)}\in\C\{t_2\}$
such that $(C,0)\cong\{(y_2,t_2)\in(\C^2,0)\,|\, 
y_2^3-\sum_{i=0}^2g_i^{(0)}t_2^i=0\}$. The system of partial
differential equations
\begin{eqnarray*}
\paa_3\begin{pmatrix}g_1\\g_0\end{pmatrix}
= \begin{pmatrix}
2g_{02}+g_{12}t_2+2g_1\\
g_{02}t_2+3g_0+\frac{2}{3}g_1g_{12}
\end{pmatrix}.
\end{eqnarray*}
is obtained from \eqref{4.45} by inserting
$(g_2,h_2,h_1,h_0)=(0,1,t_2,-\frac{2}{3}g_1)$.
By the theorem of Cauchy-Kovalevski (cited in the proof of
Corollary \ref{t4.7}), it has a unique solution with
initial values $(g_1,g_0)|_{t_3=0}=(g_1^{(0)},g_0^{(0)})$.
By construction, 
$(g_2,g_1,g_0,h_2,h_1,h_0)=(0,g_1,g_0,1,t_2,-\frac{2}{3}g_1)$
solve \eqref{4.45}. We obtain a germ of an $F$-manifold 
with $g_2=0$. It is isomorphic to a germ in Theorem \ref{t7.1}.
But now Lemma \ref{t4.6} gives
\begin{eqnarray}\label{7.26}
H_{Z_3}(Z_2)=
Z_2\cdot[2g_{22}h_2+(3y_2+g_2)h_{22}+3h_{12}]=Z_2\cdot 3.
\end{eqnarray}
The plane curve germs $(C(t_3^0),0):=(Z_2|_{t_3=t_3^0})^{-1}(0)$
are isomorphic for all $t_3^0$, but the function germs 
$Z_2|_{t_3=t_3^0}$ are not right equivalent for different
$t_3^0$, because $Z_2|_{t_3=0}$ is not quasihomogeneous
(up to coordinate changes) and because of \eqref{7.26}.
\end{remarks}

\begin{remark}\label{t7.4}
In \cite{BT14}, 3-dimensional Frobenius manifolds
with Euler fields $E=t_1\paa_1+\frac{1}{2}t_2\paa_2$
were constructed which enrich the following three
$F$-manifolds with Euler fields:
\begin{list}{}{}
\item[(i)]
The $F$-manifold $M=\C^3$ in Theorem \ref{t5.4} (a) with 
$T_tM\cong Q^{(2)}$ for all $t\in M$ and with 
Euler field $E$ as in \eqref{5.7} with 
$\varepsilon_2=\frac{1}{2}$ and $\varepsilon_{3,0}=0$. 
\item[(ii)]
The $F$-manifold $M=\C^3$ in Theorem \ref{t5.6} for $p=2$ 
(so the first one in the series) with $T_0M\cong Q^{(2)}$
and $T_tM\cong Q^{(3)}$ for generic $t\in M$ and with 
Euler field as in \eqref{5.10} with $\varepsilon_{3,0}=0$. 
\item[(iii)]
The $F$-manifold $\www M=\C^2\times\H$ which is the 
universal covering of the $F$-manifold 
$\C^2\times(\C-\{0;1\})$ in Corollary \ref{t7.2} (a)(i)
(so the one with $p=q=2$ in Theorem \ref{t7.1} (a)) 
with the Euler field as above (which is here
unique up to adding a multiple of $e$). 
\end{list}
Natural questions are now which other $F$-manifolds 
in this paper can be enriched to Frobenius manifolds
or flat $F$-manifolds, 
and with which Euler fields, and in how many ways. 
\end{remark}


\end{document}